\documentclass[conference]{IEEEtran}
\IEEEoverridecommandlockouts
\usepackage{cite}
\usepackage{amsmath,amssymb,amsfonts}
\usepackage{algorithmic}
\usepackage{graphicx}
\usepackage{textcomp}
\usepackage{xcolor}
\def\BibTeX{{\rm B\kern-.05em{\sc i\kern-.025em b}\kern-.08em
 T\kern-.1667em\lower.7ex\hbox{E}\kern-.125emX}}

\usepackage{hyperref}
\hypersetup{
 bookmarks=true,  
 pdfstartview={FitH}, 
 pdfauthor={Author}, 
 colorlinks=true, 
 linkcolor=blue,  
 citecolor=red, 
 urlcolor=cyan  
}

 \usepackage[colorinlistoftodos,prependcaption,textsize=tiny]{todonotes}

\usepackage{comment}
 
\begin{document}

\title{On the Global Convergence of Continuous–Time Stochastic Heavy–Ball Method for Nonconvex Optimization
\thanks{XZ acknowledges the Hong Kong GRF support 11305318.}
}

\author{
\IEEEauthorblockN{ Wenqing Hu}
\IEEEauthorblockA{\textit{Department of Mathematics and Statistics} \\
\textit{Missouri University of Science and Technology}\\
Rolla, USA \\
\texttt{huwen@mst.edu}}

\and

\IEEEauthorblockN{ Chris Junchi Li}
\IEEEauthorblockA{\textit{Tencent AI Lab} \\
\textit{Tecent Technology}\\
Shenzhen, China \\
\texttt{junchi.li.duke@gmail.com}}

%

\and

\IEEEauthorblockN{Xiang Zhou}
\IEEEauthorblockA{\textit{School of Data Science and Department of Mathematics, College of Science} \\
\textit{City University of Hong Kong}\\
Hong Kong, China\\
\texttt{xizhou@cityu.edu.hk}
}
}

\newtheorem{m-Theorem}{Meta-Theorem}
\newtheorem{theorem}{Theorem}[section]
\newtheorem{definition}{Definition}[section]
\newtheorem{proposition}[theorem]{Proposition}
\newtheorem{remark}[theorem]{Remark}
\newtheorem{lemma}[theorem]{Lemma}
\newtheorem{corollary}[theorem]{Corollary}
\newtheorem{example}{Example}[section]
\newtheorem{condition}{Condition}[section]
\newtheorem{assumption}{Assumption}

\providecommand{\1}{\mathbf 1}

\newcommand{\eps}{\varepsilon}
\newcommand{\abs}[1]{\left\vert#1\right\vert}
\newcommand{\set}[1]{\left\{#1\right\}}

\newcommand{\wt}[1]{\widetilde{#1}}

\maketitle

\begin{abstract}
We study the convergence behavior of the stochastic heavy-ball method with a small stepsize. Under a change of time scale, we approximate the discrete method by a stochastic differential equation that models small random perturbations of a coupled system of nonlinear oscillators. We rigorously show that the perturbed system converges to a local minimum in a logarithmic time. This indicates that for the diffusion process that approximates the stochastic heavy-ball method, it takes (up to a logarithmic factor) only a linear time of the square root of the inverse stepsize to escape from all saddle points. This results may suggest a fast convergence of its discrete-time counterpart. Our theoretical results are validated by numerical experiments.
\end{abstract}

\begin{IEEEkeywords}
heavy–ball method, dissipative nonlinear oscillator, saddle point, 
non-convex optimization, small random perturbations of Hamiltonian systems.\end{IEEEkeywords}

\section{Introduction.}

Our motivation in this work comes from the smooth unconstrained optimization problem

\begin{equation}\label{Intro:Eq:OptimizationProblem}
\min\limits_{x\in \mathbb{R}^d} f(x) \ .
\end{equation}
This problem can be solved by optimization methods that admits second--order differential equation approximations.
These methods have been demonstrated acceleration towards convergence. 
As an example, Nesterov's accelerated gradient method is a classical scheme of such type that has been used numerously in optimization.
The original method can be ``surrogated" by a limiting ODE of the form (see \cite{su2016differential})
\begin{equation}\label{Intro:Eq:ODENestrovOriginal}
\ddot{X}(t)+\dfrac{3}{t}\dot{X}(t)+\nabla_X f(X(t))=0 \ , \ X(0)\in \mathbb{R}^d \ .
\end{equation}
Here $\nabla_X f(x)$ is the gradient of the function $f$ with respect to the $X$--variable.

However, in many practices of statistical machine learning and optimization (see, e.g., \cite{[HintonMomentumICML2013]} and \cite{[AWDA],[DBSDA],[De-biasing]}), the momentum variable does not necessarily require a time--decaying factor. 
Such methods are within the range of ``heavy ball methods" date back to 1964 \cite{[Polyak1964]}, which leads to the consideration of the following version of the heavy ball method (e.g. \cite{[HintonMomentumICML2013]})
\begin{equation}\label{Intro:Eq:HeavyBallGradientVersion}
\left\{\begin{array}{ll}
x_k & =x_{k-1}+\eps v_k \ ;
\\
v_k & =(1-\alpha\eps)v_{k-1}-\eps \nabla_X f\left(x_{k-1}\right) \ .
\end{array}\right.
\end{equation}
Here under our scaling, $\alpha>0$ is the friction constant, $\eps>0$ is the learning rate which is assumed to be small, and $\mu=1-\alpha\eps \in [0,1]$ is the momentum coefficient.
It is straightforward to show that after time--rescaling $t\rightarrow t/\eps $, the family of the discrete-time processes $(x_{\lfloor t/\eps \rfloor}, v_{\lfloor t/\eps \rfloor})$ converges as $\eps \rightarrow 0$ to the solution $(X(t), V(t))$ to the system of differential equations:
\begin{equation}
\label{Intro:Eq:NonlinearOscillatorHighDimensionsHamiltonianDissipative}
\left\{\begin{array}{ll}
\dot{X}(t) & =V(t) \ ;
\\
\dot{V}(t) & = -\alpha V(t)-\nabla_X f(X(t)) \ .
\end{array}\right.
\end{equation}
Such approximation can characterized in its \textit{Ansatz} form (see e.g.~\cite[Theorem 2]{su2016differential}): for any fixed $T>0$, 
$$
 \lim\limits_{\eps\rightarrow 0}\max\limits_{0\leq t\leq T}
\left(|x_{\lfloor t/\eps \rfloor}-X(t)| +|v_{\lfloor t/\eps \rfloor}-V(t)| \right)=0
,
$$
where $| \cdot |$ is the standard Euclidean norms in $\mathbb{R}^d$.
Equation \eqref{Intro:Eq:NonlinearOscillatorHighDimensionsHamiltonianDissipative} is a Hamiltonian system with a constant friction $\alpha$.
It can be written in compact form as a second-order differential equation 
\begin{equation}\label{Intro:Eq:NonlinearOscillatorHighDimensionsDissipative}
\ddot{X}(t)+\alpha \dot{X}(t)+\nabla_X f(X(t))=0
 \ .
\end{equation}

In this work, we mainly focus on the stochastic version of \eqref{Intro:Eq:NonlinearOscillatorHighDimensionsHamiltonianDissipative} or \eqref{Intro:Eq:NonlinearOscillatorHighDimensionsDissipative} as the dynamical equation to approach local minimizers of the objective function $f(x)$.
One main reason to consider the stochastic scheme for \eqref{Intro:Eq:NonlinearOscillatorHighDimensionsHamiltonianDissipative} is to help the escape from saddle points when the objective function $f(x)$ is {\it non--convex} in optimization practice.
In this case, the deterministic process $X(t)$ in \eqref{Intro:Eq:NonlinearOscillatorHighDimensionsHamiltonianDissipative} can be trapped at local maximum points or saddle points of $f$, and therefore the question arises whether or not one can add a noise to help the escape from these unfavorable critical points. 
To this end, consider the following noisy scheme of \eqref{Intro:Eq:HeavyBallGradientVersion} which we call the \textit{stochastic heavy--ball method}:
\begin{equation}\label{Intro:Eq:HeavyBallNoisyGradientVersionNonDegenerate}
\left\{\begin{array}{ll}
x_k & =x_{k-1}+\eps \widetilde{v}_k \ ;
\\
v_k & =(1-\alpha\eps )v_{k-1}-\eps \widetilde{\nabla}_Xf\left(x_{k-1}\right) \ .
\end{array}\right.
\end{equation}
Here $\widetilde{v}_k:=v_k+\text{unbiased noise}$, and the noisy gradient $\widetilde{\nabla}_X f :=\nabla_X f+ \text{unbiased noise}$.
Let us pick the noise in the $x$--variable iteration to be independent of the noise chosen in the $v$--variable iteration. 
Following the general scheme in \cite{hu2019diffusion}, \cite{li2017stochastic,li2019stochastic}, \cite{BORKAR}, \cite{BENVENISTE-METIVIER-PRIOURET}, \cite{KUSHNER-YIN}, this leads to a system of stochastic differential equations as follows
\begin{equation}\label{Intro:Eq:WeakLimitHeavyBallNoisyNonDegenerate}
\left\{\begin{array}{ll}
dx(t) &=\eps v(t)dt+\eps \sigma_1(x(t), v(t))dW^1_t
 ;
\\
dv(t) &=-\eps(\alpha v(t) + \nabla_X f(x(t)))dt
\\
&\qquad\qquad\qquad +\eps \sigma_2(x(t), v(t))dW^2_t
.
\end{array}\right.
\end{equation}
Here  $dt$ shall be seen as 1, so it makes sense that we have both $\sqrt{s}$ and $dt$. This is a time-rescaled version of the approximating SDE seen in other literatures \cite{li2017stochastic,li2019stochastic}.
Here $W_t^1$ and $W_t^2$ are two independent $d$--dimensional Brownian motions, and $\sigma_1(x(t),v(t))$, $\sigma_2(x(t), v(t))$ originates from the noise covariance regarding to the two processes $\widetilde{v}_t$ and $\widetilde{\nabla}_X f(x)$ satisfying the non-degenerate condition that $a_i(x,v)\equiv \sigma_i (x,v)\sigma_i^T(x,v)$ are uniformly positive-definite for $i=1,2$. 
Equation \eqref{Intro:Eq:HeavyBallNoisyGradientVersionNonDegenerate} can be regarded as the time-discrertized scheme of \eqref{Intro:Eq:WeakLimitHeavyBallNoisyNonDegenerate} with a unit-time stepsize.

In this work, we shall study the stochastic continuous model \eqref{Intro:Eq:WeakLimitHeavyBallNoisyNonDegenerate} as heavy ball methods and shows that this converges to a local minimum in a linear time with respect to $\eps^{-1}$ up to a logarithmic factor.
It is well-known (e.g.~in the seminal work \cite{su2016differential}) that a continuous-time process serves as a surrogate of some discrete-time algorithms, and the convergence result of the continuous-time process \eqref{Intro:Eq:WeakLimitHeavyBallNoisyNonDegenerate} sheds light on the behaviors of the discrete-time counterpart.%
\footnote{%
We remarks on the gap of discrete algorithm and continuous-time model has been discussed and analyzed in \cite{li2017stochastic,hu2019diffusion}, who proved that for an $O(\eps^{-1})$ time interval the error gap between the discrete-time and continuous-time models is of $O(\eps)$.
Our work does not, however, aim to analyze such a gap.
}
Our main result can be formulated roughly as the following:
\begin{m-Theorem}
\label{m}
[Convergence of the diffusion limit of the stochastic heavy--ball method]
Under some mild non-degenerate condition on covariances,
as $\eps \rightarrow 0$, the process $x(t)$ in \eqref{Intro:Eq:WeakLimitHeavyBallNoisyNonDegenerate} converges to local minimizers of the objective function $f(x)$ after the time 
 $\lesssim C {\eps^{-1}} \ln(\eps^{-1})$ on average.
Here the constant $C>0$ depends on the function $f$ and the friction constant $\alpha$.
\end{m-Theorem}

To our best knowledge, our Meta Theorem \ref{m} is the first among continuous-time methods for stochastic heavy-ball method for non-convex optimization, and we compare it with another line of works discussed for discrete-time non-convex stochastic algorithms.
For instance, \cite{[ONeill-Wright],jin2017escape,jin2017accelerated} discussed the case of the so-called perturbed GD and AGD, where the gradient oracle is \textit{deterministic}.
\cite{daneshmand2018escaping} discussed fast escaping of saddle points via stochastic gradient oracles. 
Along with these existing analysis, our results offer a further hint towards the fact that the stochastic noise and acceleration are sufficiently helpful in escaping certain families of saddle points.
In comparison, our work attempts to handle the "stochastic" setting in continuous time, with stochastic gradient oracle accessible to the algorithm.
For instance, in the PGD analysis of \cite{jin2017escape} it takes (for stepsize $\eta$) $\eta^{-1}\gamma_1^{-1}$ steps to escape one saddle point, and PAGD in \cite{jin2017accelerated} uses $\eps^{-1} \sqrt{\gamma_1^{-1}}$ time to escape one saddle point, where both stepsizes are picked as $O(1)$.
Here, $\eps$ is the upper bound of the Euclidean norm of the gradient for a critical point, and $-\gamma_1$ upper bounds the least Hessian eigenvalue at the saddle point 
(see Definition \ref{Def:StrictSaddleProperty} and  \cite{ge2015escaping,sun2015nonconvex}).
Both  match the result in the stochastic setting as in this work, but our choice of stepsize needs to be much smaller in order to cope with the stochasticity. 
As the readers will see later, our analysis works with the imposed ``strong'' saddle assumption which is a quantitative characterization of the Morse function.
Removal of the strong saddle property, though, may be possible using more (nontrivial) technicalities which is left for future work.

The paper is organized as follows. In Section \ref{Sec:RandomPerturbation} we establish the connection between our continuous dynamics and the randomly perturbed dissipative nonlinear oscillators. 
In Section \ref{sec:1D} we offer a heuristic argument based on the simple calculation for the special case where $f$ is quadratic. 
Section \ref{Appendix:HamiltonianFormulation} reviews the basics of dissipative Hamiltonian system.
Section \ref{Sec:Convergence} is the main body where we provide the convergence proof by first considering the exit behavior of the randomly perturbed process near one specific saddle point and then considering global convergence in the case when there is a chain of saddle points.
In Section \ref{Sec:Numerics} we provide numerical results that validate our theory. Finally, Section \ref{Sec:Conclusion} is dedicated to some further discussions.

\section{Randomly perturbed dissipative nonlinear oscillators}\label{Sec:RandomPerturbation}

Our main idea is to reduce the problem to 
the study on the exit problem of the randomly perturbed 
Hamiltonian system \cite[Chapter 9]{[FWbookNew]}.
To illustrate the connection, consider the process $(x(t), v(t))$ defined in \eqref{Intro:Eq:WeakLimitHeavyBallNoisyNonDegenerate}.
We apply the following time rescaling technique: \begin{equation}
 \label{eqn:scaling}
 X^\eps(t):=x(t/\eps), \mbox{ and} ~ V^\eps(t):= v(t/\eps ),
 \end{equation} then
 \eqref{Intro:Eq:WeakLimitHeavyBallNoisyNonDegenerate} is equivalently 
 transformed to the following system
\begin{equation}\label{Intro:Eq:NonlinearOscillatorHighDimensionsHamiltonianRandomPerturbationEpsNonDegenerate}
\left\{\begin{array}{ll}
dX^\varepsilon(t) & =V^\varepsilon(t)dt+\sqrt{\varepsilon}\sigma_1(X^\varepsilon(t), V^\varepsilon(t))dW^1_t \ ;
\\
dV^\varepsilon(t) & =(-\alpha V^\varepsilon(t)-\nabla_X f(X^\varepsilon(t)))dt\\
&\qquad\qquad\quad +\sqrt{\varepsilon}\sigma_2(X^\varepsilon(t), V^\varepsilon(t))dW^2_t \ .
\end{array}\right.
\end{equation}

Equation \eqref{Intro:Eq:NonlinearOscillatorHighDimensionsHamiltonianRandomPerturbationEpsNonDegenerate} is a small random perturbation of \eqref{Intro:Eq:NonlinearOscillatorHighDimensionsHamiltonianDissipative}, where $\sqrt{\eps}$ represents the small intensity of the noise.
This equation is the main focus on our analytic studies in this paper.
 Since the intensity of the perturbation is small, the trajectory $(X^\varepsilon(t), V^\varepsilon(t))$ of \eqref{Intro:Eq:NonlinearOscillatorHighDimensionsHamiltonianRandomPerturbationEpsNonDegenerate} will be close to $(X(t), V(t))$
of the deterministic system \eqref{Intro:Eq:NonlinearOscillatorHighDimensionsHamiltonianDissipative} 
with the same initial condition in any finite--time interval
$[0,T]$, so that we have for any $\delta>0$ and $T>0$,
\begin{equation}\label{Intro:Eq:ConvergenceInProbabilityPerturbedSystemFiniteTime}
\lim\limits_{\varepsilon\rightarrow 0}\mathbf{P}\left(\max\limits_{0\leq t \leq T}
 (|X^\varepsilon(t)-X(t)|+|V^\varepsilon(t)-V(t)|)>\delta\right)=0
\end{equation}
under the same initial $X^\varepsilon(0)=X(0)$ and $V^\varepsilon(0)=V(0)$.

This estimation is an approximation in finite time of the perturbed process to the unperturbed one. In the
long run, e.g. with a logarithmic time scale as shown in \eqref{Thm:Kifer1981Strengthened:Eq:ExpectedExitTimeStrengthened}, the perturbed process may
escape saddle points while the unperturbed process may not be able to do so.

Before we carry out the detailed analysis for the stochastic system \eqref{Intro:Eq:NonlinearOscillatorHighDimensionsHamiltonianRandomPerturbationEpsNonDegenerate},
we first discuss some basic facts on the unperturbed system \eqref{Intro:Eq:NonlinearOscillatorHighDimensionsHamiltonianDissipative},
which is a Hamiltonian system with friction (see Appendix \ref{Appendix:HamiltonianFormulation} for more on the Hamiltonian formulation of the problem).
The Hamiltonian associated with \eqref{Intro:Eq:NonlinearOscillatorHighDimensionsHamiltonianDissipative}
is \begin{equation}\label{Intro:Eq:Hamiltonian}
{H}(X,V)=\dfrac{1}{2}V^2+f(X) \ .
\end{equation}
For any $\alpha>0$, \eqref{Intro:Eq:NonlinearOscillatorHighDimensionsHamiltonianDissipative}
 is a dissipative system in the sense that its trajectory $(X(t), V(t))$ satisfies for any $t>0$
\begin{equation}\label{Intro:Eq:DecayHamiltonianAlongDeterministicFlow}
{H}(X(t), V(t))-{H}(X(0), V(0))=-\alpha \int_0^t (V(s))^2 ds \ .
\end{equation}
From \eqref{Intro:Eq:DecayHamiltonianAlongDeterministicFlow}, we see that the
Hamiltonian function ${H}(X(t), V(t))$ is strictly decaying, unless $V(s)\equiv 0, \forall s\in (0,t)$. 
One can show that critical points
of \eqref{Intro:Eq:NonlinearOscillatorHighDimensionsHamiltonianDissipative} all lie on the $X$--axis
(or $X$--plane, and we refer this to the $X$--axis throughout the text that follows: It simply means that $V=0$).
These critical points
have $(X,V)$ coordinates $(X,0)$ in which $X$ is a critical point of the function $f$.
If $X$ is a local minimum point of the function $f(X)$,
then $(X,0)$ is a local minimum point of the Hamiltonian ${H}$;
If $X$ is a local maximum or saddle point of the function $f(X)$,
then $(X,0)$ is a saddle point of the Hamiltonian ${H}$ (see Lemma \ref{Lm:Fact1} in Section \ref{Appendix:HamiltonianFormulation}). 
 
From the above reasoning, we see the deterministic system \eqref{Intro:Eq:NonlinearOscillatorHighDimensionsHamiltonianDissipative}
to approaches a critical point of $H(X,V)$ and the trajectory $(X(t), V(t))$ may be trapped there.
 Notice that there are two types of critical points
of $H(X,V)$: saddle point or local minimum point. Due to instability of the flow near the saddle point, one can expect that random
perturbations in \eqref{Intro:Eq:NonlinearOscillatorHighDimensionsHamiltonianRandomPerturbationEpsNonDegenerate} help the process leave the saddle point
after wandering in its neighborhood for sufficiently long time.
This is exactly the reason why 
we shall devote lots of efforts below in this paper
for the randomly perturbed system \eqref{Intro:Eq:NonlinearOscillatorHighDimensionsHamiltonianRandomPerturbationEpsNonDegenerate}.
We shall see that the perturbed system $(X^\varepsilon(t), V^\varepsilon(t))$ converges only to local minimum points of $H(X,V)$ after sufficiently long time and the expected value of this convergence 
time is roughly bounded by $ O(\ln(\varepsilon^{-1}))$ for the sufficient small $\varepsilon$. See Theorem \ref{Thm:MajorResultHamiltonianForm} for details.
It is worth pointing out that when the objective function $f$ is bounded from both above and below,
then 
our conclusion holds for arbitrary initial condition.

\section{Heuristic derivation for the Quadratic Objective function}
\label{sec:1D}
To demonstrate the key ideas of our proof in Section \ref{Sec:Convergence},
we present the heuristic argument in the simple quadratic case:
 $f(x) =\frac12 x^T \Lambda x$ with $\Lambda=\mbox{diag}\{\lambda_1,\ldots, \lambda_d\}$,
 where $\lambda_i \neq 0$ for all $i$.
 $f$ now has only one saddle point $X_*=0$.
The calculation for this quadratic case is fundamentally important
since the proof of the general case is based on the local linearization assumption 
near the saddle point provided by the Hartman-Grobman Theorem
(see Assumption \ref{Assumption:Linearization} below).
Furthermore, we assume $\sigma_1=0$ to simply the calcuation. This means that 
there is no noisy term in the $v$-iteration for \eqref{Intro:Eq:HeavyBallNoisyGradientVersionNonDegenerate}. 
With these assumptions, \eqref{Intro:Eq:NonlinearOscillatorHighDimensionsHamiltonianRandomPerturbationEpsNonDegenerate}
can be reduced to 
\begin{equation}\label{eqn:quad1} 
\left\{\begin{array}{ll}
dX^\varepsilon(t) & =V^\varepsilon(t)dt \ ;
\\
dV^\varepsilon(t) & =(-\alpha V^\varepsilon(t)- \Lambda X^\varepsilon(t))dt +\sqrt{\varepsilon} S dW(t) \ .
\end{array}\right.
\end{equation}
where $SS^T=\mathbf{E}
(\wt{\nabla}_X f(X_*) [\wt{\nabla}_X f(X_*) ]^T )$
and $W_t$ is a standard $d$-dimensional Brownian motion.
Note \eqref{eqn:quad1} is a linear system. Its second-order form is 
\begin{equation}
\label{eqn:2ndXeps}
\frac{d^2 X^\eps}{ dt^2} + \alpha \frac{d X^\eps}{dt} + \Lambda X^\eps
=\sqrt{\eps} S\frac{dW}{dt}.
\end{equation}
This indicates the following rescaling in space: \begin{equation}
\wt{X} (t) := X^\eps (t) /\sqrt{\eps}, \mbox{ and } 
\wt{V} (t):=V^\eps (t)/\sqrt{\eps},
\end{equation}
then $(\wt{X} , \wt{V} )$ 
satisfies the following stochastic differential
equation independent of $\eps$:
\begin{equation}\label{eqn:linHODE} 
\left\{\begin{array}{ll}
d\wt{X}(t) & =\wt{V}(t)dt \ ;
\\
d\wt{V}(t) & =(-\alpha \wt{V}(t)- \Lambda \wt{X}(t))dt 
 + S dW(t) \ .
\end{array}\right.
\end{equation}
\eqref{eqn:2ndXeps} is then equivalent to 
$$
\frac{d^2 \wt{X}}{ dt^2} + \alpha \frac{d \wt{X} } {dt} + \Lambda 
\wt{X}
= S\frac{dW}{dt}.
$$

Write $\wt{X}=(\wt{X}_1,\ldots, \wt{X}_d)$ component-wisely.
Without loss of generality, we consider 
the first component $\wt{X}_1(t)$ of \eqref{eqn:linHODE} 
and let $e_1$ be the elementary basis vector $(1,0,\ldots,0)$
in $\mathbb{R}^d$.
We now have a linear scalar-valued stochastic
dissipative linear oscillator
for $\wt{X}_1(t)=e_1^T \wt{X}(t)$:
\begin{equation}\label{1dSDE}
\frac{d ^2 \wt{X}_1}{d t^2}
+
\alpha \frac{d \wt{X}_1}{dt} 
+\lambda_1
 \wt{X}_1 
=
e_1^T S 
\frac{d W}{d t}
=
\sigma
\frac{d \wt{W}_1}{d t},
\end{equation}
where $\sigma=\sqrt{e_1^T SS^T e_1}>0$,
and $\wt{W}_1=\sigma^{-1} e_1^T S W$
is the standard one dimensional Brownian motion.
The behavior of the solution to \eqref{1dSDE}
is mainly determined by 
the characteristic equation for a scalar $\mu$ satisfying 
$$\mu^2+\alpha \mu + \lambda_1=0.$$
Now we assume that $\lambda_1<0$ since we mainly consider the 
saddle point. Then the two roots 
\begin{equation}\label{mu}
\mu^\pm = \frac{-\alpha\pm \sqrt{\alpha^2 - 4\lambda_1}}{2}
\end{equation}
are both real and satisfies $\mu^-<0<\mu^+$.
In this case, the general solution to \eqref{1dSDE}
has the explicit form:
\begin{equation}
\begin{split}
\wt{X} _1(t)
=&
e^{ \mu^+ t }
\left(
C_1 + \sigma \int_0^t \frac{ e^{-\mu^+ s }}{\mu^+ - \mu^-} 
d\wt{W}_1(s)
\right)
\\
&
+
e^{ \mu^- t }
\left(
C_2 - \sigma \int_0^t \frac{ e^{-\mu^- s }}{\mu^+ - \mu^-} 
d\wt{W}_1(s)
\right)
\\
& \sim \ e^{ \mu^+ t }.
\end{split}
\end{equation}
One can also verity that the variance 
$\mbox{Var}(\wt{X}_1(t)) \sim e^{ 2\mu^+ t }
$ for large $t$.

Now recall that the original process $x(t)$ in \eqref{Intro:Eq:WeakLimitHeavyBallNoisyNonDegenerate}
is linked to $\wt{X}(t)$ here by the following scaling both in time and space
$x(t)=X^\eps(\eps t)=\sqrt{\eps} \wt{X}(\eps t)$.
So, $x_1(t)=\sqrt{\eps} \wt{X}_1(\eps t)\sim \sqrt{\eps} e^{\mu^+ \eps t} $.
To exit a neighbour around the saddle point $X_*=0$,
it will takes $x_1(t)$ the time roughly about $ \frac12 \frac{1}{\mu^+} \eps^{-1} \ln (\eps^{-1})$, which is 
consistent with 
the rate given in Meta-Theorem \ref{m}.

We hope this simple analysis for the escape behavior near the saddle point
can shed light on the insight for our rigorous analysis below. 
The main result is Theorem \ref{Thm:ExitOneSaddlePoint} 
where there are two key ingredients in the proof. The first one is the linearization analysis for the $2d$ dimensional 
dissipative Hamiltonian system
for all possible values of the eigenvalues $\lambda_i$, 
which is a generalization of the above eigenvalue analysis \eqref{mu} for \eqref{1dSDE}.
Refer to Proposition \ref{prop:LinearAlgebraTechnical} in Appendix \ref{Appendix:StructureHamiltonianFlowOneSpecificSaddle}.
The second ingredient is to apply the main result of \cite{[Kifer1981]} which 
was originally for the first order dynamical system in our current setup. 
Our method is to consider the Hamiltonian flow in the $(x,v)$ phase space 
$\mathbb{R}^{2d}$.

\section{Hamiltonian Formation of the Dissipative Nonlinear Oscillator}
\label{Appendix:HamiltonianFormulation}

In this section, we provide standard Hamiltonian formulations for our randomly perturbed dissipative oscillator system.
Throughout the text, $\nabla_X$ or $\nabla_V$ will denote gradient with respect to $X$ or $V$ variable, respectively, and $\nabla^2_X$, $\nabla^2_V$,
etc. are denoted similarly;
$\nabla$  denotes gradient with respect to $(X,V)$ variable, and $\nabla^2$ are denoted similarly; if we use $\dfrac{\partial}{\partial X}$
or $\dfrac{\partial}{\partial V}$, then it means the corresponding gradients with respect to $X$ or $V$ variable, and $\dfrac{\partial^2}{\partial X^2}$
or $\dfrac{\partial^2}{\partial V^2}$, etc. are defined similarly. Standard Euclidean norms in $\mathbb{R}^d$ are defined either by $|\bullet|_{\mathbb{R}^d}$
or $|\bullet|$.

We consider a Hamiltonian system
\begin{equation}\label{Eq:HeavyBallStrongLimitHamiltonianWithoutFriction}
\left\{\begin{array}{ll}
dX(t) & =V(t)dt \ , \ X(0)\in \mathbb{R}^d \ ;
\\
dV(t) & =-\nabla_X f(X(t))dt \ , \ V(0)\in \mathbb{R}^d \ .
\end{array}\right.
\end{equation}
We add friction to \eqref{Eq:HeavyBallStrongLimitHamiltonianWithoutFriction}, and we get
\begin{equation}\label{Eq:HeavyBallStrongLimitHamiltonian}
\left\{\begin{array}{ll}
dX(t) & =V(t)dt \ , \ X(0)\in \mathbb{R}^d \ ;
\\
dV(t) & =(-\alpha V(t)-\nabla_X f(X(t)))dt \ , \ V(0)\in \mathbb{R}^d \ .
\end{array}\right.
\end{equation}

We will denote the flow map of \eqref{Eq:HeavyBallStrongLimitHamiltonian} to be $S^t$, so that
$(X(t), V(t))=S^t(X(0), V(0))$. Let us define the Hamiltonian 
\begin{equation}\label{Eq:Hamiltonian}
H(X,V)=\dfrac{1}{2}V^2+f(X) \ ,
\end{equation}
so that
$$\dfrac{\partial H}{\partial V}=V \text{ and }
\dfrac{\partial H}{\partial X}=\nabla_X f(X) \ .$$
%
%
%

Let $(X_O, V_O)$ be a critical point of $H$. The above implies that $V_0=0$ and $X_0$
is a critical point of $f(X)$.
From \eqref{Eq:Hamiltonian} we have a formal expansion
\begin{equation}\label{Eq:FormalExpansionHamiltonianCriticalPoint}
\begin{split}
H(X,V)=&H(X_0, 0)+\dfrac{1}{2}V^2+(X-X_0)^T \nabla_X^2 f(X) (X-X_0)
\\
&~\qquad + \text{higher order terms in } X .
\end{split}
\end{equation}
If $X_0$ is a local minimum point of $f(X)$, then \eqref{Eq:FormalExpansionHamiltonianCriticalPoint} tells us that $(X_0, 0)$
is a local minimum point of $H(X,V)$. If $X_0$ is a local maximum or saddle point of $f(X)$, then
\eqref{Eq:FormalExpansionHamiltonianCriticalPoint} tells us that $(X_0, 0)$ is a saddle point of $f(X)$. In particular, by our
strong saddle property assumption for the potential function $f$
(see Definition \ref{Def:StrongSaddleProperty} in Section \ref{Sec:Convergence}), 
it is easy to see that the Hamiltonian function $H(X,V)$
also has the strong saddle property. 
The summary of these discussions
 is  presented in Lemma \ref{Lm:Fact1}.
 
Define the column vector $Y(t)=\begin{pmatrix} X(t) \\ V(t)\end{pmatrix}\in \mathbb{R}^{2d} $ and the skew gradient
\begin{equation}\label{Eq:SkewGradientHamiltonian}
\nabla^\perp H(X,V)=\begin{pmatrix}
\dfrac{\partial H}{\partial V}
\\
-\dfrac{\partial H}{\partial X}
\end{pmatrix} \ .
\end{equation}
Then the system
\eqref{Eq:HeavyBallStrongLimitHamiltonian} can be written in a standard Hamiltonian form with friction

\begin{equation}\label{Eq:HeavyBallStrongLimitHamiltonianStandard}
\dfrac{dY}{dt}=\nabla^{\perp} H(Y(t))+b(Y(t)) \ , Y(0)=Y_0\in \mathbb{R}^d\times \mathbb{R}^d \ .
\end{equation}

Here $b(X,V)=(0, -\alpha V)^T$ is the friction vector field. Notice that $\text{div} b=-\alpha<0$ corresponds
to the classical friction case. Furthermore, $b(X,V)$ is a gradient field in the sense that we can write
$b(X,V)=-\nabla B(X,V)$, in which $B(X,V)=\dfrac{1}{2}\alpha V^2$.

Let us introduce the standard symplectic matrix

\begin{equation}\label{Eq:StandardSymplecticMatrix}
J=\begin{pmatrix}
0 & I_d
\\
-I_d & 0\end{pmatrix} \ .
\end{equation}
The standard symplectic matrix $J$ has the property that $J^2=-I_{2d}$, $J^T=-J=J^{-1}$.

The Hamiltonian vector field $\nabla^\perp H(Y)$ can be written as $\nabla^\perp H(Y)=J\nabla H(Y)$. The local behavior of this vector field near its critical point is characterized by the ``skew Hessian matrix''

\begin{equation}\label{Eq:SkewHessianMatrix}
\nabla [J \nabla H(X,V)]
=
\begin{pmatrix}
\dfrac{\partial^2 H}{\partial X\partial V} & \dfrac{\partial^2 H}{\partial V^2}
\\
-\dfrac{\partial^2 H}{\partial X^2} & -\dfrac{\partial^2 H}{\partial X \partial V}\end{pmatrix}
=
\begin{pmatrix} 0 & I_d \\ -\nabla_X^2 f(X) & 0 \end{pmatrix} \ ,
\end{equation}
in the sense that we have the expansion around a critical point $(X_O, V_O)$ of $\nabla^\perp H$:
\begin{equation}\label{Eq:ExpansionHamiltonianFieldCriticalPoint}
\begin{split}
\nabla^\perp H (X,V)&=\nabla J \nabla H(X_O,V_O)\begin{pmatrix} X-X_O \\ V-V_O\end{pmatrix}
\\
+&\psi(X-X_O, V-V_O)
(|X-X_O|_{\mathbb{R}^d}^2+|V-V_O|_{\mathbb{R}^d}^2) \ ,
\end{split} 
\end{equation}
where $\psi(X,V)$ is some bounded smooth vector--function in the variables $(X,V)$.

Let us also calculate the Hessian matrix of the Hamiltonian function $H$ as

\begin{equation}\label{Eq:HessianMatrix}
\nabla^2 H(X,V)
=
\begin{pmatrix}
\dfrac{\partial^2 H}{\partial X^2} & \dfrac{\partial^2 H}{\partial X\partial V}
\\
\dfrac{\partial^2 H}{\partial X \partial V} & \dfrac{\partial^2 H}{\partial V^2} \end{pmatrix}
=
\begin{pmatrix} \nabla_X^2 f(X) & 0 \\ 0 & I_d \end{pmatrix} \ .
\end{equation}

From \eqref{Eq:SkewHessianMatrix} and \eqref{Eq:HessianMatrix}, taking into account that $\nabla_X^2 f$ is a symmetric matrix,
 we see that we have the relation

\begin{equation}\label{Eq:RelationSkewHessianMatrixANDHessianMatrix}
J\nabla^2 H(X,V)= \nabla[ J \nabla H(X,V) ] .
\end{equation}

Let us denote the matrix

\begin{equation}\label{Eq:I0}
I^0=
\begin{pmatrix}
0 & 0
\\
0 & I_d
\end{pmatrix}
 \ .
\end{equation}
Then we can write the friction term $b(X,V)$ as
\begin{equation}\label{Eq:FrictionMatrixRepresentation}
b(Y)=-\alpha I^0 Y
\end{equation}
where $Y=\begin{pmatrix}
X \\ V
\end{pmatrix}$.
%
%
%
%
%
%
Then the linearized behavior of the dissipative Hamiltonian oscillator (21) is determined by the Jacobi matrix
\begin{equation}\label{A}
A
=
\nabla [J\nabla H(X,V)] - \alpha I^0
=
\begin{pmatrix}
0, & I_d
\\
-\nabla_X^2 f(X), & -\alpha I_d
\end{pmatrix}.
\end{equation}

\section{Rigorous Proof of Global convergence in general case}
\label{Sec:Convergence}

Before we continue with the analysis, we will make some additional structural assumptions on the function $f(x)$.

\subsection{Strong saddle property}

\begin{definition}\label{Def:Morse}
A smooth function $f:\mathbb{R}^d\to \mathbb{R}$ is a {\it Morse function} if it has all critical points being {\it non--degenerate}, i.e.,
the Hessian $\nabla^2 f(x_O)$ at any critical point $x_O$ is non--degenrate. This implies that all eigenvalues
$\lambda_1\leq ...\leq \lambda_d$ at $x_O$ are nonzero.
\end{definition}

Morse functions admit local quadratic re--parametrization at each critical point,
which is the content of the so--called {\it Morse Lemma} \cite[Lemma 2.2]{[Milnor1963]}.
To ensure that the perturbation helps the process $X(t)$ escape from saddle points,
we introduce the following ``strict saddle property"
(compare with \cite{ge2015escaping}, \cite{sun2015nonconvex}) as follows.

\begin{definition}[strict saddle property]\label{Def:StrictSaddleProperty}
Given fixed $\gamma_1>0$ and $\gamma_2>0$, we say a Morse function $f$ defined on $\mathbb{R}^d$ satisfies the
``strict saddle property" if each point $x\in \mathbb{R}^d$ belongs to one of the following three cases:
\begin{enumerate}
\item[(i)] $|\nabla f(x)|> \gamma_2>0$ \ ;
\item[(ii)] $|\nabla f(x)|\leq \gamma_2$ and $\lambda_{\min}(\nabla^2 f(x))\leq -\gamma_1 <0$ \ ;
\item[(iii)] $|\nabla f(x)|\leq \gamma_2$ and $\lambda_{\min}(\nabla^2 f(x))\geq \gamma_1>0$.
\end{enumerate}
Here $\lambda_{\min}(\nabla^2 f(x))$ is the minimal eigenvalue of the Hessian matrix $\nabla^2 f(x)$ at point $x$.
\end{definition}

We will call a saddle point $x\in \mathbb{R}^d$ of the function $f$ a ``strict saddle" if Definition \ref{Def:StrictSaddleProperty} (ii)
 holds at $x$. 
 Thus a Morse function $f$
that satisfies the strict saddle property has all its saddle points that are strict saddle points.
Note that the strict saddle property above only focus on the minimal eigenvalues 
but does not imply anything on the 
degeneracy of the eigenvalues.

For the sake of proof, we need to assume that all eigenvalues of the Hessian $\nabla^2 f(x)$ at critical points are 
{\it uniformly} bounded away from
$0$ in absolutely value
(which is stronger than the non-degeneracy requirement in the definition of Morse function). This leads to our new notion of ``strong saddle property" as follows.

\begin{definition}[strong saddle property]\label{Def:StrongSaddleProperty}
 We say the Morse function $f$ satisfies the
``strong saddle property" if it satisfies the strict saddle property defined above
{ and} there exits a constant $\gamma_3>0$ such that 
for any critical point $x_O$ ( i.e., $\nabla f(x_O)=0$), 
all eigenvalues $\lambda_i (x_O) $, $i=1,2,...,d$, of the Hessian $\nabla^2 f(x_O)$ satisfy 
 $ |\lambda_i (x_O)| \geq \gamma_3>0$.
\end{definition}

We will call a saddle point $x\in \mathbb{R}^d$ of the function $f$ a ``strong saddle" if
Definition \ref{Def:StrongSaddleProperty} holds at $x$. Note that strong saddle points are the strict saddle points where the absolute values of all
eigenvalues of the Hessian are bounded away from 0 by a positive constant $\gamma_3$.
Throughout this paper we will work under Definition \ref{Def:StrongSaddleProperty} for the objective function $f$.

\subsection{Exit behavior near one specific saddle point}\label{Sec:Convergence:OneSaddle}

To carry out the analysis like in Section \ref{sec:1D} and to apply the main theories in \cite{[Kifer1981]},
we write the system \eqref{Intro:Eq:NonlinearOscillatorHighDimensionsHamiltonianRandomPerturbationEpsNonDegenerate} in 
the standard Hamiltonian form (see Section \ref{Appendix:HamiltonianFormulation})
\begin{equation}\label{Eq:WeakLimitHamiltonianStandardRandomPerturbationEps:Main}
dY^\varepsilon_t=[\nabla^\perp H(Y^\varepsilon_t)+b(Y^\varepsilon_t)]dt+\sqrt{\eps}\Sigma(Y^\varepsilon_t)dW_t \ . \end{equation}
where $$Y^\varepsilon_t: =(X^\varepsilon(t), V^\varepsilon(t))^T\in \mathbb{R}^{2d}$$ 
is seen as the column vector in $\mathbb{R}^{2d}$.

We study in this subsection the exit behavior of \eqref{Eq:WeakLimitHamiltonianStandardRandomPerturbationEps:Main} near one specific saddle point. Our method here follows that of \cite{[Kifer1981]}, \cite{[MikamiExitUnstable]}, \cite{[DayExitSaddle]}, \cite{[BakhtinHeteroclinic]},
\cite{hu2017fast}, 
among many other literature dedicated to this topic. The main difference here is
 that we are working with the Hamiltonian dynamics.

In the case when $\alpha>0$ is small as the parameter $\varepsilon>0$ goes to zero,
the general program of dealing with a randomly perturbed Hamiltonian system such as \eqref{Eq:WeakLimitHamiltonianStandardRandomPerturbationEps:Main} is considered in \cite{[Freidlin-WeberExit]},
\cite{[Brin-Freidlin]}, \cite{[Freidlin-WeberAnnProb]}, \cite{[Freidlin-WeberExit]}, \cite{[Freidlin-WeberSPA]}, \cite[Chapter 9]{[FWbookNew]}. In our case, the friction term $b(X,V)=(0, -\alpha V)$ has a magnitude that is moderate when compared to the Hamiltonian flow $\nabla^\perp H$,
and it does not go to zero together with the small parameter $\varepsilon > 0$. In this case, one has to carefully analyze the behavior of the process \eqref{Eq:WeakLimitHamiltonianStandardRandomPerturbationEps:Main} near the saddle point.

Recall than the deterministic version of \eqref{Eq:WeakLimitHamiltonianStandardRandomPerturbationEps:Main}
is the Hamiltonian flow with friction defined in \eqref{Intro:Eq:NonlinearOscillatorHighDimensionsHamiltonianDissipative}.
We study the critical points and the linear part of \eqref{Intro:Eq:NonlinearOscillatorHighDimensionsHamiltonianDissipative}.
First of all, we have the following lemma.
\begin{lemma}\label{Lm:Fact1}
Let $\alpha> 0$. All critical points of the system \eqref{Intro:Eq:NonlinearOscillatorHighDimensionsHamiltonianDissipative} lie on the $X$--axis.
These critical points are all saddle points of the Hamiltonian flow $\nabla^\perp H$ and 
have $(X,V)$ coordinates $(X,0)$ where $X$ is a critical point of $f(X)$.
If $X$ is a local minimum point of the function $f(X)$,
then $(X,0)$ is a local minimum point of the Hamiltonian $H$;
If $X$ is a local maximum or saddle point of the function $f(X)$,
then $(X,0)$ is a saddle point of the Hamiltonian $H$. All saddle points of the Hamiltonian function $H(X,V)$
are strong saddle points in Definition \ref{Def:StrongSaddleProperty}.
\end{lemma}

\

By Lemma \ref{Lm:Fact1}, let $O=(X_O,0)\in \mathbb{R}^{2d}$ be a saddle point of the Hamiltonian.
then $X_O$ is a local maximum or saddle point of the function $f$.
The linear part of the flow \eqref{Intro:Eq:NonlinearOscillatorHighDimensionsHamiltonianDissipative} 
near the saddle point $O$ is given by the $2d\times 2d$ Jacobi matrix $A$ defined in \eqref{A}.
Since $\nabla^2_X f(X_O)$ is a symmetric matrix, we can find an orthonormal basis $\xi_1,...,\xi_d$ (viewed as column vectors) in $\mathbb{R}^d$ such that $\nabla^2_X f(O) \xi_i=\lambda_i \xi_i$. $X_O$ is a local maximum or saddle point of $f(x)$, and thus
we see that without loss of generality we can assume that
\begin{equation}\label{Eq:SpectrumOfgradf}
\lambda_1\leq \lambda_2\leq...\leq \lambda_k<0<\lambda_{k+1}\leq...\leq \lambda_{d}
\end{equation}
for some $1\leq k\leq d$.
 $k$ is the number of negative eigenvalues,
the index of the saddle point of $X_O$.

Proposition \ref{prop:LinearAlgebraTechnical} in Appendix \ref{Appendix:StructureHamiltonianFlowOneSpecificSaddle} 
is an important tool bridging the eigenvalues of 
the matrix $A$ for the dissipative Hamiltonian flow \eqref{Intro:Eq:NonlinearOscillatorHighDimensionsHamiltonianDissipative} with $(\lambda_i, \xi_i)$ the eignpairs of the Hessian $\nabla^2 f(X_O)$.
This proposition is the substantial development of our one dimensional calculation in Section \ref{sec:1D}
and the cornerstone of our main results.

To state our main result 
Theorem \ref{Thm:ExitOneSaddlePoint},
we present the following specifications.
Let $G$ be a connected open neighborhood in $\mathbb{R}^{2d}$ of the saddle point $O$
with the smooth boundary $\partial G$, such that $O$ is the only critical point of
the Hamiltonian flow $\nabla^\perp H$ inside $G$. 
 Let the process $Y_t^\varepsilon=(X^\varepsilon(t), V^\varepsilon(t))^T$ defined in
\eqref{Intro:Eq:NonlinearOscillatorHighDimensionsHamiltonianRandomPerturbationEpsNonDegenerate} start from initial condition $Y_0^\varepsilon=(x,v) \in G$. Let
\begin{equation}\label{Eq:HittingTimeNeighborhoodSaddle}
\tau^\varepsilon_{(x,v)}=\inf\{t>0: Y^\eps_0=(x,v) \ , \ Y^\varepsilon_t\in \partial G\} \ .
\end{equation}

Denote the flow map of \eqref{Intro:Eq:NonlinearOscillatorHighDimensionsHamiltonianDissipative}
be defined as $S^t$. Introduce the decomposition as in \cite{[Kifer1981]}:
$$G\cup \partial G=0\cup A_1\cup A_2\cup A_3 \ ,$$
where $A_1$ is a set of points $(x,v)\in G\cup \partial G$ such that if $(x,v)\in A_1$ then $S^u(x,v)\in G$
for $u>s$ and $S^u (x,v)\not \in G\cup \partial G$ if $u\leq s$ for some $s=s(x,v)\leq 0$ and $S^t (x,v)\rightarrow 0$ as $t\rightarrow \infty$;
$A_2$ is a set of points $(x,v)\in G\cup \partial G$ such that if $(x,v)\in A_2$ then $S^u (x,v)\in G$ for $u<s$
and $S^u (x,v)\not\in G\cup \partial G$ if $u>s$ for some $s=s(x,v)\geq 0$ and $S^t (x,v)\rightarrow 0$ as $t\rightarrow -\infty$;
$A_3$ is a set of points $x\in G\cup \partial G$ such that if $(x,v)\in A_3$ then $S^u (x,v)\in G$ provided $s_1<u<s_2$
and $S^u (x,v)\not\in G\cup \partial G$ if either $u>s_2$ or $u<s_1$ for some $s_1=s_1(x,v)\leq 0$ and $s_2=s_2(x,v)\geq 0$.
In other words, $A_1$ is the set of initial points which will enter $G$ from outside,
 $A_2$ is the set of initial points which will exit $G$ from inside,
 and $A_3$ is the set of initial points which will cross (i.e., enter and then exit) $G$.
If $(x,v)\in A_2\cup A_3$, then $S^t (x,v)$ leaves $G$ after some time, so that there is a finite
\begin{equation}\label{Eq:DeterministicFiniteExitTime}
t(x,v)=\inf\{t>0 : S^t (x,v)\in \partial G\} \ .
\end{equation}

To specify the exit distribution on $\partial G$, 
we need some further technical assumptions as in \cite{[Kifer1981]}.
Assume there exists some $1\leq k^\circ \leq k$ such that $\lambda_1=\lambda_2=...=\lambda_{k^\circ}$. 
So $k^\circ$ is the multiplicity of the lowest eigenvalues. 
Denote by $\gamma_{\text{max}}$ the eigenspace of $A$ in \eqref{A} which corresponds to the eigenvalues
$\mu_1^+,...,\mu_{k^\circ}^+$. Then as in \cite{[Kifer1981]}, there exists a $k^\circ$--dimensional sub--manifold $W_{\text{max}}$
tangent to $\gamma_{\text{max}}$ at the saddle point $O$ and is invariant with respect to $S^t$. We see that $Q_{\text{max}}=W_{\text{max}}\cap \partial G$
is not empty. 

Now we state our first theorem. Define 
\begin{equation}\label{Eq:mu0:Main}
\mu_0=\dfrac{-\alpha+\sqrt{\alpha^2-4\lambda_1}}{2}>0 .
\end{equation}
which plays the similar role to \eqref{mu}.
\begin{theorem}\label{Thm:ExitOneSaddlePoint}
(a) (Exit time) For each fixed $(x,v)\in G$ we have
\begin{equation}\label{Thm:ExitOneSaddlePoint:Eq:ExitTimeAsymptotics}
\limsup\limits_{\varepsilon\rightarrow 0}\dfrac{\mathbf{E} \tau^\varepsilon_{(x,v)}}{\ln (\varepsilon^{-1})}\leq \dfrac{1}{2\mu_0} .
\end{equation}

(b) (Exit distribution) If $(x,v)\in (O \cup A_1)\backslash \partial G$ then for any open set $Q$ of $\partial G$ such that
$Q \supset Q_\text{max}$ we have
\begin{equation}\label{Thm:ExitOneSaddlePoint:Eq:ExitDistributionStableManifold}
\lim\limits_{\varepsilon\downarrow 0}\mathbf{P}(Y_{\tau_{(x,v)}^\varepsilon}^\varepsilon\in Q)=1 \ .
\end{equation}
If $(x,v)\in A_2\cup A_3$ then for any Borel measurable set $Q$ of $\partial G$ we have
\begin{equation}\label{Thm:ExitOneSaddlePoint:Eq:ExitDistributionUnstableManifold}
\lim\limits_{\varepsilon\downarrow 0}\mathbf{P}(Y_{\tau_{(x,v)}^\varepsilon}^\varepsilon \in Q)=\1(S^{t(x,v)}(x,v) \in Q) \ .
\end{equation}

\end{theorem}

\

The proof is the direct application of Theorems 2.1--2.3 of \cite{[Kifer1981]}
in view of Proposition \ref{prop:LinearAlgebraTechnical} in Appendix \ref{Appendix:StructureHamiltonianFlowOneSpecificSaddle} and we skip the proof.

Theorem \ref{Thm:ExitOneSaddlePoint} characterizes the asymptotic exit time for the exit from one specific saddle point
and the bound is in the sense of super limit. 
In the next, we extend this theorem to a uniform bound in Theorem \ref{Thm:Kifer1981Strengthened}.
We would like to point out here, that improving Theorem \ref{Thm:ExitOneSaddlePoint} to Theorem \ref{Thm:Kifer1981Strengthened}
requires essential exploitation of the Linearization Assumption (Assumption \ref{Assumption:Linearization} below).
By the classical Hartman--Grobman Theorem (see \cite[\S 13]{[Arnold]}), for any strong saddle point $O$
that we consider, there exists an open neighborhood $G$ of $O$, and a $\mathbf{C}^{(0)}$ homeomorphism mapping $h: G\rightarrow \mathbb{R}^n$,
such that the dissipative Hamiltonian flow given by \eqref{Intro:Eq:NonlinearOscillatorHighDimensionsHamiltonianDissipative} is mapped by $h$ into a linear flow.
The homeomorphism $h$ is called
a (linear) conjugacy mapping. To make our argument work, we will have to put a stronger assumption.

\begin{assumption}[Linearization Assumption]\label{Assumption:Linearization}
The homeomorphism $h$ provided by the Hartman-Grobman Theorem
can be taken to be $\mathbf{C}^{(2)}$.
\end{assumption}

It is known that a sufficient condition for the
validity of the $\mathbf{C}^{(2)}$ Linearization Assumption is the so called non--resonance condition (see, for example,
the Sternberg linearization Theorem \cite[Theorem 6.6.6]{[KatokHasselblatt]}).

Let $U$ be an open neighborhood of the saddle point $O$ such that
$\text{dist}(U\cup \partial U, \partial G)>0$. 
 The following theorem holds under Assumption \ref{Assumption:Linearization}.
\begin{theorem}\label{Thm:Kifer1981Strengthened}
(a) (Exit time)
For any $r>0$, there exist some $\varepsilon_0>0$ so that for all $(x,v)\in U\cup \partial U$ and all $0<\varepsilon<\varepsilon_0$ we have
\begin{equation}\label{Thm:Kifer1981Strengthened:Eq:ExpectedExitTimeStrengthened}
\dfrac{\mathbf{E}\tau_{(x,v)}^\varepsilon}{\ln(\varepsilon^{-1})}\leq \dfrac{1}{2\mu_0}+r \ .
\end{equation}
Here the stopping time $\tau_{(x,v)}^\varepsilon$ is defined as in \eqref{Eq:HittingTimeNeighborhoodSaddle}, and $\mu_0$
is defined as in \eqref{Eq:mu0:Main}.

(b) (Exit distribution)
For any small $\kappa>0$ and any $\rho>0$, there exist some $\varepsilon_0>0$ so that for all $x\in U\cup \partial U$ and all $0<\varepsilon<\varepsilon_0$
we have
\begin{equation}\label{Thm:Kifer1981Strengthened:Eq:ExitDistributionStrenghtened}
\mathbf{P}(Y_{\tau_{(x,v)}^\varepsilon}^\varepsilon\in Q^\kappa)\geq 1-\rho \ .
\end{equation}
Here
$Q^\kappa=\{(x,v)\in \partial G \ , \ \text{dist}((x,v), \partial G_{U\cup \partial U\rightarrow \text{out}})<\kappa\}$
with
\begin{equation*}
\begin{split}
\partial G_{U\cup \partial U \rightarrow {\text{out}}}=
&
~Q_{\text{max}}\cup \{S^{t(x,v)}(x,v) \text{ for some }
\\
&(x,v)\in U \cup \partial U \text{ with finite } t(x,v)\} \ .
\end{split}
\end{equation*}
\end{theorem}

\

The {proof} of this Theorem can be found in Appendix \ref{Appendix:ProofOfTheoremKifer1981Strengthened}.

\subsection{Chain of saddle points and the global convergence}\label{Sec:Convergence:ManySaddles}

Upon exit from a neighborhood of a saddle point, the process $Y^\varepsilon_t$ in \eqref{Eq:WeakLimitHamiltonianStandardRandomPerturbationEps:Main}
may further involve multiple saddle points before reach a local minimum.
We will analyze the case where there is 
 a chain of saddle points in this part. We would like to point out that there might be an infinite number of saddle points in the landscape. Nevertheless, our results apply to the
convergence time when the continuous dynamics passes $k$ consecutive saddle points among all of these saddle
points, so that we are analyzing the time during which the dynamics visits the $k$ consecutive saddle points.

Assume that there are $k$--consecutive strong saddle points
$O_1,...,O_k$ of the Hamiltonian $H(X,V)$ such that $H(O_1)>H(O_2)>...>H(O_k)$. Let $O^*$ be a local minimum point of
$H(X,V)$ such that $H(O_k)>H(O^*)$ and there are no other critical points $o$ such that $H(O_k)>H(o)>H(O^*)$.

Notice that any critical point of the Hamiltonian flow $\nabla^\perp H$ has $V$--component equal to zero. Furthermore,
the saddle points $O_1,...,O_k$ are of the form $O_j=(x^{O_j},0)$, $j=1,2,...,k$ where $x^{O_j}$ is a local maximum point or saddle point
(under Definition \ref{Def:StrongSaddleProperty}) of the potential function $f(x)$. In the same way we have
$O^*=(x^*, 0)$ where $x^*$ is a local minimum point of the potential function $f(x)$.

Select a small number $e>0$.
For the perturbed process $Y_t^\varepsilon$ in \eqref{Eq:WeakLimitHamiltonianStandardRandomPerturbationEps:Main},
 define the stopping time

\begin{equation}\label{Eq:HittingTimeChainOfSaddle}
\mathcal{T}_{(X,V)}^{H,\varepsilon}=\inf
\left\{t\geq 0: Y_0^\varepsilon=(X,V), H(Y_t^\varepsilon)<H^*+e\right\} \ .
\end{equation}
where $H^*:=H((x^*,0))= f(x^*)$ is the minimum of $H$
(as well as of $f$) at the local minimizer $(x^*,0)$.
Similarly, we define
\begin{equation}\label{Eq:HittingTimeChainOfSaddleFunctionf}
\mathcal{T}_{(X,V)}^{f,\varepsilon}=\inf\{t\geq 0: X_0^\varepsilon=X \ , \ f(X_t^\varepsilon)<f(x^*)+e\}.
\end{equation}
It is clear that
 $\mathcal{T}_{(X,V)}^{H,\varepsilon} \geq \mathcal{T}_{(X,V)}^{f,\varepsilon}$
 since $H(Y_t^\eps)=f(X_t^\eps)+ \frac12 |V_t^\eps|^2$ and $H^*=f(x^*)$.
In this section we aim at characterizing the asymptotic upper bound as $\varepsilon\rightarrow 0$ of the times $\mathcal{T}_{(X,V)}^{H,\varepsilon}$
(or $\mathcal{T}_{(X,V)}^{f,\varepsilon}$).

For each saddle point $O_i$, we consider a nested pair of open neighborhoods $U_i \subsetneq G_i$ containing $O_i$.
Let us first pick
$k$ disjoint neighbours $G_i$, $i=1,2,...,k$ in such a way that $O_i\in G_i$ is the only critical point inside $G_i$, and for any $i\neq j$ we have $G_i \cap G_j=\emptyset$.
Then select $U_i\subset G_i$
such that $\text{dist}(\partial U_i, \partial G_i)>0$ for all $i=1,2,...,k$. For each saddle point $O_i$,
let us denote by $\gamma_{i, \text{max}}$ the eigenspace of the Jacobi matrix$A(O_i)$
which corresponds to the
eigenvalues with the largest positive real parts (see \eqref{A} and Theorem \ref{Thm:ExitOneSaddlePoint}) .
For each open neighborhood $G_i$,
as in \cite{[Kifer1981]}, we can construct a submanifold $W_{i, \text{max}}$ that is tangent
to $\gamma_{i, \text{max}}$ at $O_i$ and is invariant with respect to $S^t$.
Let $Q_{i, \text{max}}=W_{i, \text{max}}\cap \partial G_i$. By the classification of points in $G_i\cup \partial G_i$
as $O_i$, $A_1$, $A_2$,
$A_3$ as in \cite{[Kifer1981]} presented right above Theorem \ref{Thm:ExitOneSaddlePoint}, we see that for any point
$(x,v)\in U_i\cup \partial U_i$,
either there exist some finite $t(x,v)$ such that $S^{t(x,v)}(x,v)\in \partial G_i$, or $S^t (x,v)\rightarrow O_i$ as $t\rightarrow\infty$.
Let
\begin{equation*}
\begin{split}
\partial G_{i, U_i\cup \partial U_i \rightarrow {\text{out}}}=
&~ Q_{i, \text{max}}\cup \{S^{t(x,v)}(x,v) \text{ for some }
\\
& (x,v)\in U_i \cup \partial U_i \text{ with finite } t(x,v)\} \ .
\end{split}
\end{equation*}
 
For the small $e>0$, define 
\begin{equation*}
\begin{split}
\mathcal{U}=\{(x,v)\in \mathbb{R}^{2d}: H(x,v)<H(O^*)+e 
\\
 ~\qquad \text{ for some local minimum point } O^*\} \ .
\end{split}
\end{equation*}
The set
$\mathcal{U}=\bigcup\limits_j \mathcal{U}_j$
where each $\mathcal{U}_j$ is an open neighborhood of one of the local minimum points of the Hamiltonian $H(X,V)$.
We can set $e>0$ to be so small that each of $\mathcal{U}_i$ and $\mathcal{U}_j$ are mutually disjoint, and they are also disjoint with any of
the $G_i$'s.

Let us pick the neighborhoods $G_i$ sufficiently small such that starting from any point $(x,v)$ on $G_i \cup \partial G_i$, the deterministic flow $S^t(x,v)$ of \eqref{Intro:Eq:NonlinearOscillatorHighDimensionsHamiltonianDissipative}
will never return to $G_i\cup \partial G_i$. This can be achieved by the friction term $b(X,V)$. Let $U=\bigcup\limits_{i=1}^k U_i$ and $G=\bigcup\limits_{i=1}^k G_i$. Let

\begin{equation}\label{Eq:InfSkewGradOutsideU}
K=\inf\left\{|\nabla^\perp H(x,v)|: (x,v)\not\in U\cup \mathcal{U}\right\} \ .
\end{equation}

From our construction above we see that $K>0$ and it is independent of $\varepsilon$. 
With the initial $Y_0^\varepsilon=(X_0,V_0)^T$,
we define the sequence of stopping times
\begin{equation*}
0=\sigma_0^\eps\leq \tau_1^\eps\leq \sigma_1^\eps\leq \tau_2^\eps\leq \sigma_2^\eps\leq ...
\end{equation*}
where \begin{equation*}
\tau_j^\eps:=\inf\left\{t>\sigma_{j-1}^\eps: Y_t^\varepsilon\in \partial U\cup \partial \mathcal{U}\right\} \ ,
\end{equation*}
and
\begin{equation*}
\sigma_j^\eps:=\inf\left\{t>\tau_{j}^\eps: Y_t^\varepsilon\in \partial G \cup \partial \mathcal{U}\right\} \ .
\end{equation*}

Starting from any initial condition $Y^\eps_0=(X_0, V_0)$ outside of $G\cup \mathcal{U}$, such that 
there exists constants $H_1, H_2$
we have 
\begin{equation}\label{Hb}
H_1<H(X_0, V_0)<H_2,
\end{equation}
the process $Y_t^\varepsilon$ travels for time $\tau_1^\eps$ before it enters
$\overline{U\cup \mathcal{U}}$. 
We can bound the expected time $\tau_1^\eps$ 
in the following proposition (Lemma 3.3 in \cite{hu2017fast}).

\begin{proposition}\label{prop:LyapunovFunctionFiniteTravelTime}
There exists some $\varepsilon_0>0$ uniformly for all $Y^\eps_0=(X_0, V_0)$ with $H_1<H(X_0,V_0)<H_2$,
such that for all $0<\varepsilon<\varepsilon_0$, \begin{equation}\label{Eq:TravelTimeFinite}
\mathbf{E}_x \tau_1^\eps \leq C \ .
\end{equation}
for a finite number $C>0$ independent of $\varepsilon$.
\end{proposition}
\

Finally by using the arguments of \cite[Section 3]{hu2017fast}, we then have
the theorem below as a generalization of 
Theorem \ref{Thm:ExitOneSaddlePoint}.

\begin{theorem}\label{Thm:MajorResultHamiltonianForm}
Consider the process $Y_t^\varepsilon$ defined as in \eqref{Eq:WeakLimitHamiltonianStandardRandomPerturbationEps:Main}
with the initial conditions bounded by \eqref{Hb}.
Assume the diffusion matrix $\Sigma(Y)\Sigma^T(Y)$ is uniformly positive definite for all choices of $Y$.
Then we have

(i) For any small $\rho>0$, with probability at least $1-\rho$, the process $Y_t^\varepsilon$ in
\eqref{Eq:WeakLimitHamiltonianStandardRandomPerturbationEps:Main} converges
to the local minimum point $O^*$ for sufficiently small $\varepsilon$ after passing through all $k$ saddle points $O_1$, ..., $O_k$;

(ii) As $\varepsilon\downarrow 0$, conditioned on the above convergence of $Y_t^\varepsilon$ to $O^*$,
we have
\begin{equation}\label{Thm:MajorResultHamiltonianForm:Eq:TotalExitTimeAsymptotic}
\limsup\limits_{\varepsilon \rightarrow 0}\dfrac{\mathbf{E} \mathcal{T}_{(X,V)}^{H,\varepsilon}}{\ln (\varepsilon^{-1})}\leq
\dfrac{k}{4\gamma_1}\left(\sqrt{\alpha^2+4\gamma_1}+\alpha\right) \ .
\end{equation}
Here the stopping time $\mathcal{T}_{(X,V)}^{H,\varepsilon}$ is defined in \eqref{Eq:HittingTimeChainOfSaddle}, and $\gamma_1$ is in Definition \ref{Def:StrictSaddleProperty}.
\end{theorem}

The {proof} of Proposition \ref{prop:LyapunovFunctionFiniteTravelTime}
 is in Appendix \ref{Appendix:ProofPropositionLyapunovFunctionFiniteTravelTime} 
and the 
 {proof} of Theorem \ref{Thm:MajorResultHamiltonianForm} is in Appendix \ref{Appendix:ProofTheoremMajorResultHamiltonianForm}.

 \subsection{Main Theorem}
 With the preparation of the previous three subsections,
 we now are ready to state our main theorem.

 Theorem \ref{Thm:MajorResultHamiltonianForm} is formulated for the Hamiltonian function $H(X,V)=\dfrac{1}{2}V^2+f(X)$. 
 It is straightforward to reformulate it in term of $x$ variable and the objective function $f$.

\begin{corollary}\label{Cor:MajorResultOptimizationForm}
Consider the process $X_t^\varepsilon$ defined as in \eqref{Intro:Eq:NonlinearOscillatorHighDimensionsHamiltonianRandomPerturbationEpsNonDegenerate}
with the initial position and momentum bounded by \eqref{Hb}.
Let $x^*$ be the unique local minimum of
$f$ within an open neighborhood $U(x^*)$ such that $f(x^*)<f(x^{O_k})$;
refer to the beginning of Section \ref{Sec:Convergence:ManySaddles} for the saddle points $O_k$.
 Then we have

(i) For any small $\rho>0$, with probability at least $1-\rho$, the process $X_t^\varepsilon$ in
\eqref{Intro:Eq:NonlinearOscillatorHighDimensionsHamiltonianRandomPerturbationEpsNonDegenerate} converges
to the minimizer $x^*$ for sufficiently small $\varepsilon$ after passing through all $k$ saddle points $x^{O_1}$, ..., $x^{O_k}$;

(ii) As $\varepsilon\downarrow 0$, conditioned on the above convergence of $X_t^\varepsilon$ to $x^*$,
we have
\begin{equation}\label{Cor:MajorResultOptimizationForm:Eq:TotalExitTimeAsymptotic}
\limsup\limits_{\varepsilon \rightarrow 0}\dfrac{\mathbf{E} \mathcal{T}_{(X,V)}^{f,\varepsilon}}{\ln (\varepsilon^{-1})}\leq
\dfrac{k}{4\gamma_1}\left(\sqrt{\alpha^2+4\gamma_1}+\alpha\right) \ .
\end{equation}
Here $\mathcal{T}_{(X,V)}^{f,\varepsilon}$ is defined as in \eqref{Eq:HittingTimeChainOfSaddleFunctionf}.
\end{corollary}

\

Finally, we formulate the convergence result for the diffusion approximation of the stochastic heavy ball method
\eqref{Intro:Eq:WeakLimitHeavyBallNoisyNonDegenerate}. Taking into account that the process $(x(t), v(t))$
in \eqref{Intro:Eq:WeakLimitHeavyBallNoisyNonDegenerate} is related to the process $(X(t), V(t))$
in \eqref{Intro:Eq:NonlinearOscillatorHighDimensionsHamiltonianRandomPerturbationEpsNonDegenerate} via a time change 
$$(X(t), V(t))=(x(t/\eps), v(t/\eps)).$$
 we see that we have the following Theorem
 immediately from Corollary \ref{Cor:MajorResultOptimizationForm}.

\begin{theorem}\label{Thm:MajorResultOptimizationFormHeavyBall}
Consider the process $x(t)$ defined as in \eqref{Intro:Eq:WeakLimitHeavyBallNoisyNonDegenerate}.
 Let the objective function $f(x)$ satisfy the strong saddle property in Definition \ref{Def:StrongSaddleProperty}.
Let $x^*$ be the unique local minimum of
$f$ within an open neighborhood $U(x^*)$ such that $f(x^*)<f(x^{O_k})$. Then we have

(i) For any small $\rho>0$, with probability at least $1-\rho$, the process $x(t)$ in
\eqref{Intro:Eq:WeakLimitHeavyBallNoisyNonDegenerate} converges
to the minimizer $x^*$ for sufficiently small $s>0$ after passing through all $k$ saddle points $x^{O_1}$, ..., $x^{O_k}$;

(ii) Set $e>0$ small and let
\begin{equation}\label{Thm:MajorResultOptimizationFormHeavyBall:Eq:HittingTimeChainOfSaddle}
T_x=\inf\{t\geq 0: x(0)=x \ , \ f(x(t))<f(x^*)+e\} \ .
\end{equation}
Then as $s\downarrow 0$, conditioned on the above convergence of $x(t)$ to $x^*$,
we have
\begin{equation}\label{Thm:MajorResultOptimizationFormHeavyBall:Eq:TotalExitTimeAsymptotic}
\limsup\limits_{\eps \rightarrow 0}\dfrac{\mathbf{E} T_x}{{\eps^{-1}}\ln (\eps^{-1})}\leq
\dfrac{k}{4\gamma_1}\left(\sqrt{\alpha^2+4\gamma_1}+\alpha\right) \ .
\end{equation}
for any initial $x$ whose function value $f(x)$ is bounded $H_1<f(x)<H_2$. 
$\gamma_1$ is in Definition \ref{Def:StrictSaddleProperty}.

\end{theorem}

\begin{remark}
Using similar analysis proposed in this paper, which also dates back to \cite{[Kifer1981]},  the first hitting time to a neighborhood of local minimizers for continuous-time SGD (with stepsize $\eta$) is asymptotically bounded by $k \gamma_1^{-1} \eta^{-1}\log \eta^{-1}$, and our analysis for accelerated SGD (with stepsize $\eps$) reduces it to $\asymp k \gamma_1^{-0.5} \eps^{-1} \log \eps^{-1}$. 
Compared to \cite{hu2017fast}, such a result indicates that using the same stepsize, stochastic heavy-ball method also escapes from all saddles and helps the iteration to reach the local minimum point at a reduced period of time by $\gamma_1^{-0.5}$, showing its comparative advantages for saddle-point escaping when $\gamma_1$ is relatively small.
\end{remark}

\section{Numerical results}
\label{Sec:Numerics}

\begin{figure}[!tb]
\centering
\includegraphics[width=0.42\textwidth]{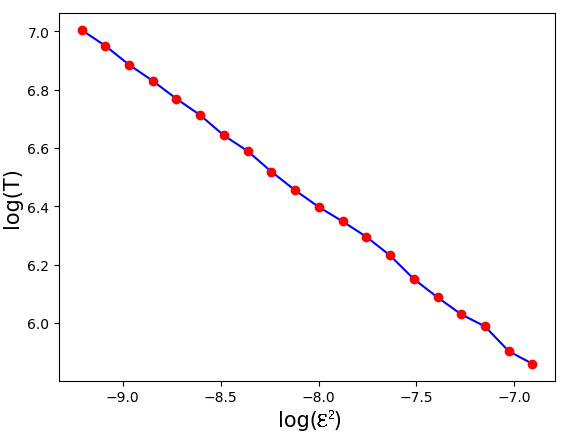}
\caption{The relationship between the learning rate $s=\eps^2$ and the stopping time $T_x$.}
\label{fig:relation_stepsize}
\end{figure}
\begin{figure}[!tb]
\centering
\includegraphics[width=0.42\textwidth]{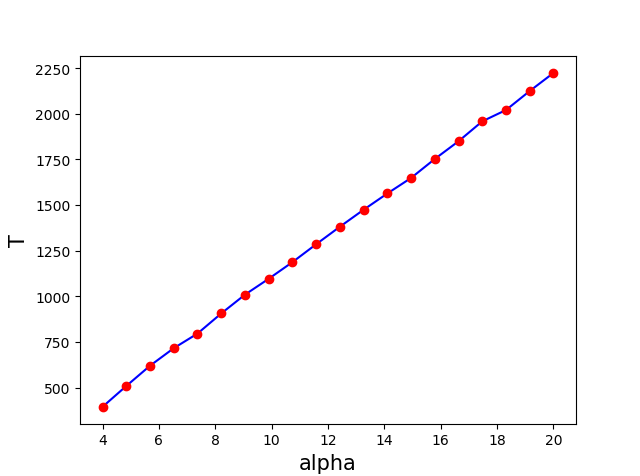}
\includegraphics[width=0.42\textwidth]{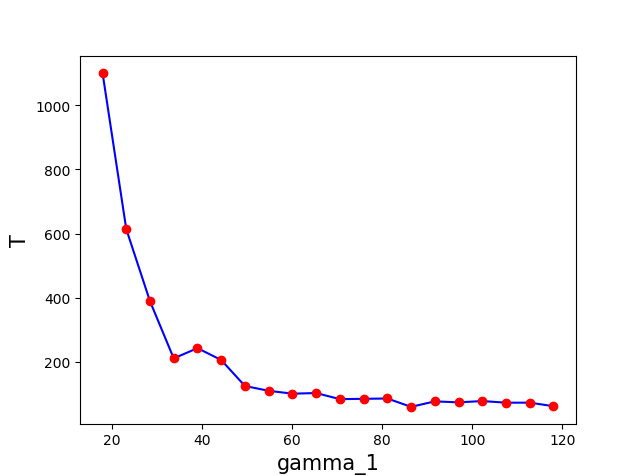}
\caption{The relationship between the parameters $\alpha$, $\gamma_1$ and the stopping time $T_x$.}
\label{fig:relation_params}
\end{figure}
\begin{figure}[!tb]
\centering
\includegraphics[width=0.4\textwidth]{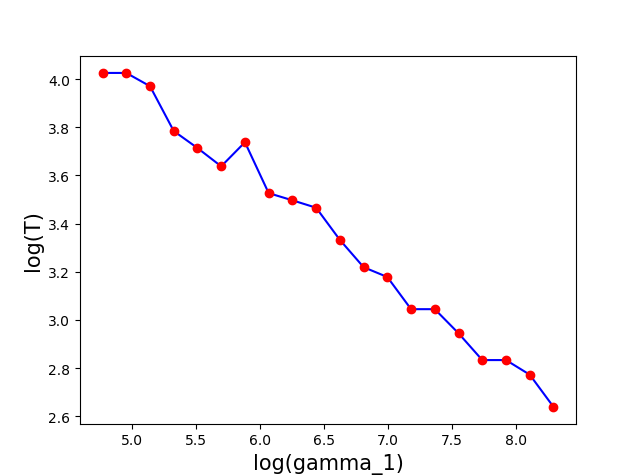}
\caption{The relationship between the parameter $\gamma_1 \equiv \alpha^2$ and the stopping time $T_x$.}
\label{fig:relation_params2}
\end{figure}

In this section, we verify the relationship between the stepsize $\eps$, the friction constant $\alpha$, the parameter $\gamma_1$ and the convergence time $T_x$, which was previously given in \eqref{Thm:MajorResultOptimizationFormHeavyBall:Eq:TotalExitTimeAsymptotic} in Theorem \ref{Thm:MajorResultOptimizationFormHeavyBall}.
We emphasize that there is a strong connection between our main result and algorithms implemented in discrete time. To be clear, we view the algorithmic scheme as a discretized SDE, and hence when the stepsize $\eps$ is sufficiently small the scheme can be viewed as numerical solution to the SDE.

We apply the stochastic heavy ball method to a cubic-regularized quadratic function of form $f(x) = x^T \Lambda x + \|x\|_2^3$. For simplicity, we consider dimension $d = 10$ and the diagonal $\Lambda$, whose eigenvalues are strictly larger than 0 or strictly smaller than 0.
Our dimension considered is 10, and the example of $\Lambda$ used in simulation code is,
$$
\Lambda = \text{diag}\{9, 7, 5, 3, 1, -1, -3, -5, -7, -\gamma_1 / 2\}
, $$
where $\gamma_1 = 18$ by default and may vary accordingly.

We can choose the parameter $\gamma_1 = - 2 \lambda_{\min} > 0$. We consider the noise in stochastic heavy ball method following unbiased random normal distributions. Convergence time $T_x$ is measured by the first time $f(x_t) - f(x^*) < C$, where $x^*$ is the local minimizer and $C = 10^{-3}$ is a fixed positive threshold. We properly initialized $x_0$ near the saddle point 0.

We first verify the relationship between stepsize $\eps$ and convergence time $T_x$, by holding parameters $\alpha$ and $\gamma_1$ fixed. Figure~\ref{fig:relation_stepsize} plots the logarithm value of stopping time $T_x$ 
with respect to logarithm of different time stepsize $s=\eps^2$. It can be observed that $\log(T_x)$ has a negative linear relationship with $\log(s)$, with slope being approximately $- 1 / 2$. This agrees with our theoretical time complexity of $T_x \asymp  \eps^{-1} \ln (\eps^{-1})$ given in \eqref{Thm:MajorResultOptimizationFormHeavyBall:Eq:TotalExitTimeAsymptotic} up to a logarithmic factor $\ln (\eps^{-1})$ neglected.

Now to investigate how the change in paramters $\alpha$ and $\gamma_1$ affects convergence time, we vary each parameter while fixing the other one and stepsize. From the theoretical time complexity, we would expect to see a positive relationship between $\alpha$ and $T_x$, and a negative relationship between $\gamma_1$ and $T_x$. This is verified by numerical experiment results shown in Figure~\ref{fig:relation_params}. As shown in the left figure, the relationship between $\alpha$ and convergence time $T_x$ appears to be linear, satisfying the asymptotic linearity implied by theoretical time complexity \eqref{Thm:MajorResultOptimizationFormHeavyBall:Eq:TotalExitTimeAsymptotic} when $\alpha >\!\!> \gamma_1$. In the right figure of $\gamma_1$ and convergence time $T_x$, we observed a close to inverse relationship, which is explained by the theoretical time complexity \eqref{Thm:MajorResultOptimizationFormHeavyBall:Eq:TotalExitTimeAsymptotic}.

To drill further into the constant factor in \eqref{Thm:MajorResultOptimizationFormHeavyBall:Eq:TotalExitTimeAsymptotic}, we vary $\alpha$ and $\gamma_1$ together by fixing $\gamma_1 \equiv \alpha^2$, and demonstrate the relationship. Then based on the theoretical time complexity \eqref{Thm:MajorResultOptimizationFormHeavyBall:Eq:TotalExitTimeAsymptotic}, the constant factor in our case is $\dfrac{k}{4\gamma_1}\left(\sqrt{\alpha^2+4\gamma_1}+\alpha\right) = \dfrac{1+\sqrt{5}}{4} \gamma_1^{-0.5}$ and we expect to see a linear relationship between $\log(\gamma_1)$ and $\log(T_x)$ with slope $- 0.5$, which is well displayed in Figure~\ref{fig:relation_params2}.

Together, Figures~\ref{fig:relation_stepsize} and~\ref{fig:relation_params} exhibit numerical results that are in accordance with \eqref{Thm:MajorResultOptimizationFormHeavyBall:Eq:TotalExitTimeAsymptotic},
and we have substantiated our main Theorem \ref{Thm:MajorResultOptimizationFormHeavyBall} from a numerical perspective.

\section{Conclusion and further discussions}\label{Sec:Conclusion}

Our work connects the behavior of the stochastic heavy-ball method with a stochastic differential equation that describes small random perturbations of a coupled system of nonlinear oscillators. By showing that the perturbed system converges to local minimizers in logarithmic time, we conclude that the continuous--stochastic heavy ball method takes (up to logarithmic factor) only a linear time of the square root of inverse stepsize to evade from all saddle points and hence it implies fast convergence of the continuous stochastic heavy--ball method.


\bibliographystyle{plain}

\newpage\appendices
\onecolumn

\section{The structure of the Hamiltonian flow near one specific saddle point}
\label{Appendix:StructureHamiltonianFlowOneSpecificSaddle}

In this Appendix we study the structure of the Hamiltonian flow near one specific saddle point. According to \eqref{Eq:ExpansionHamiltonianFieldCriticalPoint}, the linear part of the Hamiltonian vector field $\nabla^\perp H=J\nabla H$ should presumably be given by $\nabla J \nabla H(O)$ as in \eqref{Eq:SkewHessianMatrix}. By \eqref{Eq:RelationSkewHessianMatrixANDHessianMatrix}, we know that
$\nabla J \nabla H(O)=J \nabla^2 H(O)$, which is not necessarily a symmetric matrix, and in fact it will even produce complex eigenvalues that may have either or both nonzero real and imaginary parts (think of the symplectic matrix $J$ itself, the eigenvalues are $\pm i$). In this case, using \eqref{Eq:ExpansionHamiltonianFieldCriticalPoint}, we know that the linear part of
the Hamiltonian vector field $\nabla^\perp H$ with dissipation near the saddle point $O$ is given by the matrix 
$$A=JQ-\alpha I^0$$
Here $Q$ is the symmetric Hessian matrix $Q=\nabla^2 H (O)$.

Let us also note that in \cite{[BakhtinHeteroclinic]} the author has to assume that all the eigenvalues of $JQ$ are real and simple, and in our case it is natural to expect complex eigenvalues. This will be a difference between our case and the case considered in \cite{[BakhtinHeteroclinic]}, where
all eigenvalues of the matrix $A$ are real and simple.
Due to this reason, the linearization idea used in equation (8.1) of \cite{[BakhtinHeteroclinic]}
needs some further investigation, which we postpone the discussion to the next section.
However, even in the case of complex eigenvalues (with all eigenvalues having non--zero real parts),
the argument of \cite{[Kifer1981]} still works, but is only able to help us analyze the exit behavior near one specific saddle point.
To apply this argument, in our case for the system
\eqref{Eq:WeakLimitHamiltonianStandardRandomPerturbationEps:Main}, we have to consider (technically) the eigenvalues and
eigenvectors of the matrix $A=JQ-\alpha I^0$.

Let us first imagine that $I^0$ is replaced by a full rank identity matrix $I=I_{2d}$. Then all complex eigenvalues of $A=JQ-\alpha I$
differ from the ones of $JQ$ by a shift $-\alpha$ in their real parts. As a consequence, if $\alpha>0$ is very large, it is possible that none
of the complex eigenvalues of $A$ have positive real part. In this case, the exit to the boundary of a neighborhood of the saddle point $O$
takes exponentially long time due to large deviation effects (see \cite[Chapter 4]{[FWbookNew]}). Geometrically, we see that the projection of the unstable manifold of $\nabla^\perp H$
to the $X$ direction will always create an unstable component pointing from the saddle $O$. Given the half--rank matrix $I^0$, we see that the friction term
 $b(Y)=-\alpha I^0 Y$
has zero projection to the $X$--direction. Such intuitive reasoning can be improved into the following technical result in linear algebra.

%

Let the saddle point $O$ be $O=(X_0,0)$.
Since $\nabla^2_X f(X_0)$ is a symmetric matrix, we can find an orthonormal basis $\xi_1,...,\xi_d$ (viewed as column vectors) in $\mathbb{R}^d$
such that $\nabla^2_X f(O) \xi_i=\lambda_i \xi_i$. By Lemma \ref{Lm:Fact1}, $X_0$ is a local maximum or saddle point of $f(x)$, and thus
we see that without loss of generality
we can assume that
\begin{equation}\label{Eq:SpectrumOfgradf}
\lambda_1\leq \lambda_2\leq...\leq \lambda_k<0<\lambda_{k+1}\leq...\leq \lambda_{d}
\end{equation}
for some $1\leq k\leq d$.

\begin{proposition}\label{prop:LinearAlgebraTechnical}
There exist an invertible $2d\times 2d$ matrix $P$ with real or complex terms
such that for the matrix $A=JQ-\alpha I^0$ we have
$$P^{-1}AP=\text{diag}(A_1,...,A_d)$$
with each $A_i$ being a $2\times 2$ block matrix.
Moreover, the invertible matrix $P$ can be taken of the form
$$P=(\boldsymbol{u}_1^+ , \boldsymbol{u}_1^- , ... , \boldsymbol{u}_d^+, \boldsymbol{u}_d^-) \ ,$$
in which $\boldsymbol{u}_i^{\pm}$ are $2d$--dimensional real or complex vectors for $i=1,2,...,d$.

Here for some integers $l,m$ such that $1\leq k\leq l\leq m\leq d$ we have
\begin{enumerate}

\item For $i=1,...,k$, we have $\lambda_i<0$,
   $A_i=\begin{pmatrix}\mu_i^+ & 0\\ 0& \mu_i^-\end{pmatrix}$,
   $\mu_i^{\pm}=\dfrac{-\alpha \pm \sqrt{\alpha^2-4\lambda_i}}{2}$,
   $\mu_i^{\pm}$ are real, $\mu_i^+>0>\mu_i^-$,
   and $\boldsymbol{u}_i^{\pm}=\begin{pmatrix}\xi_i \\ \mu_i^{\pm}\xi_i \end{pmatrix}$ are real $2d$--dimensional vectors;

\item For $i=k+1,...,l$, we have $0<\lambda_i<\dfrac{\alpha^2}{4}$,
  $A_i=\begin{pmatrix}\mu_i^+ & 0\\ 0& \mu_i^-\end{pmatrix}$, $\mu_i^{\pm}$ are real, $0>\mu_i^+>\mu_i^-$,
  $\mu_i^{\pm}=\dfrac{-\alpha \pm \sqrt{\alpha^2-4\lambda_i}}{2}$,
  and $\boldsymbol{u}_i^{\pm}=\begin{pmatrix}\xi_i \\ \mu_i^{\pm}\xi_i \end{pmatrix}$ are real $2d$--dimensional vectors;

\item For $i=l+1,...,m$, we have $\lambda_i=\dfrac{\alpha^2}{4}$,
  $A_i=\begin{pmatrix}-\frac{\alpha}{2} & 1\\ 0& -\frac{\alpha}{2}\end{pmatrix}$, $\mu_i^+=\mu_i^-=-\dfrac{\alpha}{2}$,
  and $\boldsymbol{u}_i^{+}=\begin{pmatrix}\xi_i\\ -\frac{\alpha}{2}\xi_i\end{pmatrix}$, $\boldsymbol{u}_i^-=\textbf{a}_i$
  are two linearly independent real $2d$--dimensional vectors, in which $\textbf{a}_i$ satisfies
  $(A-\lambda_i I) \textbf{a}_i=\begin{pmatrix}\xi_i \\ -\frac{\alpha}{2}\xi_i \end{pmatrix}$;

\item For $i=m+1,...,d$, we have $0<\dfrac{\alpha^2}{4}<\lambda_i$,
  $A_i=\begin{pmatrix}\mu_i^+ & 0\\ 0& \mu_i^-\end{pmatrix}$, $\mu_i^{\pm}$ are complex and are of the form
  $\mu_i^{\pm}=-\dfrac{\alpha}{2}\pm \dfrac{\sqrt{4\lambda_i-\alpha^2}}{2}i$,
  and $\boldsymbol{u}_i^{\pm}=\begin{pmatrix}\xi_i \\ \mu_i^{\pm}\xi_i \end{pmatrix}$ are complex $2d$--dimensional vectors.
\end{enumerate}

\end{proposition}

\begin{IEEEproof}
Consider the matrix
$$A=JQ-\alpha I^0=\begin{pmatrix} 0 & I_d \\ -\nabla^2_X f(X_0) & -\alpha I_d \end{pmatrix} \ .$$
Suppose an eigenvector of the matrix $A$ has the form $\begin{pmatrix} \xi \\ v \end{pmatrix}$ with eigenvalue $\mu$ ($\mu$ may be complex).
Here $\xi$ and $v$ are two column vectors in $\mathbb{R}^d$.
Then we have

\[
\begin{split}
 A\begin{pmatrix} \xi \\ v\end{pmatrix}
=\begin{pmatrix} 0 & I_d \\ -\nabla^2_X f(X_0) & -\alpha I_d \end{pmatrix}\begin{pmatrix} \xi \\ v\end{pmatrix}
\\
=\begin{pmatrix} v \\ -\nabla_X^2 f(X_0) \xi -\alpha v\end{pmatrix}
=\mu \begin{pmatrix} \xi \\ v\end{pmatrix} \ .\end{split}
\]

This implies that $v=\mu \xi$ and $-\nabla_X^2 f(X_0) \xi=(\mu+\alpha)v=(\mu+\alpha)\mu \xi$. Therefore $\xi$ must be an eigenvector of $\nabla^2_X f(X_0)$
with eigenvalue $\lambda=(\mu+\alpha)\mu$.

Conversely, if $\xi$ is an eigenvector of $\nabla^2_X f(X_0)$ with eigenvalue $\lambda$, say
$\nabla^2_X f(X_0) \xi=\lambda \xi$, then $\begin{pmatrix}\xi \\ \mu \xi\end{pmatrix}$ is an eigenvector of $A=JQ-\alpha I^0$
with $A\begin{pmatrix}\xi \\ \mu \xi \end{pmatrix}=\mu\begin{pmatrix}\xi \\ \mu \xi \end{pmatrix}$ and the eigenvalue $\mu$ satisfies
$\mu(\mu+\alpha)=\lambda$.

From the above we see that there is a correspondence between eigenvectors/eigenvalues of $\nabla_{X}^2 f(X_0)$
and eigenvectors/eigenvalues of $A=JQ-\alpha I^0$. In fact, each eigenvector $\xi_i$ ($i=1,2,...,d$) of $\nabla_X^2 f(X_0)$ with eigenvalue $\lambda_i$
corresponds to two eigenvectors $\begin{pmatrix} \xi \\ \mu_i^{+}\xi \end{pmatrix}$ and
$\begin{pmatrix} \xi \\ \mu_i^{-}\xi \end{pmatrix}$ of $A=JQ-\alpha I^0$ with eigenvalues $\mu_i^+$ and $\mu_i^-$
\footnote{It can happen that $\mu_i^+=\mu_i^-$, and in that case the two eigenvectors may alternatively be replaced
by a two--dimensional invariant subspace. We will discuss this case later in this proof.}.

The two eigenvalues $\mu_i^{\pm}$ are the two roots of the equation

\begin{equation}\label{prop:LinearAlgebraTechnical:Eq:EquationOfMu}
\mu_i^2+\alpha \mu_i+\lambda_i=0 \ ,
\end{equation}
so that
\begin{equation}\label{prop:LinearAlgebraTechnical:Eq:EquationOfMuRoots}
\mu_i^{\pm}=\dfrac{-\alpha \pm \sqrt{\alpha^2-4\lambda_i}}{2} \ .
\end{equation}

Let us analyze the eigenvalues $\mu_i^{\pm}$ from \eqref{prop:LinearAlgebraTechnical:Eq:EquationOfMuRoots}.
Recall that we have $\lambda_1\leq \lambda_2\leq...\leq \lambda_k<0<\lambda_{k+1}\leq...\leq \lambda_{d}$. We discuss the following cases:

\begin{enumerate}
\item For $i=1,...,k$, we have $\lambda_i<0$. The two eigenvalues $\mu_i^{\pm}$ are real and $\mu_i^+>0>\mu_i^-$.

\item For $i=k+1,...,l$, we have $0<\lambda_i<\dfrac{\alpha^2}{4}$ and $\mu_i^{\pm}$ are real and $0>\mu_i^+>\mu_i^-$.

\item For $i=l+1,...,m$, we have $\lambda_i=\dfrac{\alpha^2}{4}$ and $\mu_i^+=\mu_i^-=-\dfrac{\alpha}{2}$;

\item For $i=m+1,...,d$, we have $0<\dfrac{\alpha^2}{4}<\lambda_i$, and $\mu_i^{\pm}$ are complex and are of the form
  $$\mu_i^{\pm}=-\dfrac{\alpha}{2}\pm \frac{\sqrt{4\lambda_i-\alpha^2}}{2}i \ ;$$
\end{enumerate}

In summary, the only eigenvalues of $A=JQ-\alpha I^0$ that have positive real parts
are $\mu_1^+, ..., \mu_k^+$, and all the other eigenvalues of $A=JQ-\alpha I^0$ have negative real parts.

Recall that $\lambda_1<0$ is
the negative eigenvalue with largest absolute value among all eigenvalues $\lambda_1\leq\lambda_2\leq...\leq \lambda_k<0$. If we set
$\mu_0=\dfrac{-\alpha+\sqrt{\alpha^2-4\lambda_1}}{2}$, then we have $\max\limits_{i=1,2,...,d} \text{Re} \mu_i^{\pm}=\mu_0>0$.

Notice that when $\alpha^2\neq 4\lambda_i$, the two eigenvalues $\mu_i^+\neq \mu_i^-$, and thus the eigenvectors
$$\begin{pmatrix}\xi_i \\ \mu_i^+ \xi_i\end{pmatrix} \ , \ \begin{pmatrix}\xi_i \\ \mu_i^- \xi_i\end{pmatrix}$$
are linearly independent. In fact, if we have
$$c_i^+\begin{pmatrix}\xi_i \\ \mu_i^+\xi_i\end{pmatrix}
+c_i^-\begin{pmatrix}\xi_i \\ \mu_i^-\xi_i\end{pmatrix}
=0 \ , $$
then we have
$$\left\{\begin{array}{l} c_i^++c_i^-=0 \ , \\ c_i^+\mu_i^++c_i^-\mu_i^-=0 \ ,\end{array}\right.$$
 from which we derive $c_i^+=c_i^-=0$ under $\mu_i^+\neq \mu_i^-$. In this case the two dimensional linear invariant
 subspace
$V_i=\text{span}\left\langle \begin{pmatrix}\xi_i \\ \mu_i^+\xi_i \end{pmatrix} ,
\begin{pmatrix}\xi_i \\ \mu_i^-\xi_i \end{pmatrix}\right\rangle$ splits further into two independent $1$--dimensional subspaces
$V_i=V_i^+ \oplus V_i^-$, with $V_i^{\pm}=\text{span}\left\langle \begin{pmatrix}\xi_i \\ \mu_i^{\pm}\xi_i \end{pmatrix}\right\rangle$.
The corresponding Jordan block for the invariant subspace $V_i$ is diagonal $\begin{pmatrix} \mu_i^+& 0\\ 0 & \mu_i^-\end{pmatrix}$.

If for some $i=k+1,...,d$ we have $\alpha^2=4\lambda_i$, then from the eigenvector
$\begin{pmatrix}\xi_i \\ -\frac{\alpha}{2}\xi_i \end{pmatrix}$ and the real eigenvalue $\mu_i^\pm=-\dfrac{\alpha}{2}$ of $A=JQ-\alpha I^0$
we solve the generalized eigenvalue/eigenvector equation
$$(A-\lambda_i I) \textbf{a}_i=\begin{pmatrix}\xi_i \\ -\frac{\alpha}{2}\xi_i \end{pmatrix} \ ,$$
and we obtain another vector $\textbf{a}_i\in \mathbb{R}^{2d}$ indepenent of $\begin{pmatrix}\xi_i \\ -\frac{\alpha}{2}\xi_i \end{pmatrix}$.
The two--dimensional linear subspace
$V_i=\text{span}\left\langle \begin{pmatrix}\xi_i \\ -\frac{\alpha}{2}\xi_i \end{pmatrix}, \textbf{a}_i\right\rangle$
forms an invariant subspace of the matrix $A=JQ-\alpha I^0$ for the eigenvalue $\mu_i^\pm=-\dfrac{\alpha}{2}$,
which corresponds to the Jordan block $\begin{pmatrix} -\frac{\alpha}{2} & 1\\ 0 & -\frac{\alpha}{2}\end{pmatrix}$.

Notice that if for some $i\neq j$ and $i,j\in \{k+1,...,d\}$ we have $\alpha^2=4\lambda_i=4\lambda_j$, then the eigenvectors
$\begin{pmatrix}\xi_i \\ -\frac{\alpha}{2}\xi_i \end{pmatrix}$ and $\begin{pmatrix}\xi_j \\ -\frac{\alpha}{2}\xi_j \end{pmatrix}$
are linearly independent since $\xi_i$ and $\xi_j$ are linearly independent. This implies that the invariant subspace for the matrix $A=JQ-\alpha I^0$
corresponding to
the eigenvalue $-\dfrac{\alpha}{2}$ with a possible multiplicity splits into two--dimensional invariant subspaces as described
in the previous paragraph.

Summarizing the above discussion, we see that we can find an invertible matrix of the form
\[ 
\begin{split}
P&=
(
\boldsymbol{u}_1^+ , \boldsymbol{u}_1^- , ... , \boldsymbol{u}_k^+ , \boldsymbol{u}_k^- , \boldsymbol{u}_{k+1}^+ , 
\boldsymbol{u}_{k+1}^-,
\\
 &\quad ..., 
\boldsymbol{u}_{l}^+, \boldsymbol{u}_{l}^-, \boldsymbol{u}_{l+1}^+, \boldsymbol{u}_{l+1}^-, 
...,\boldsymbol{u}_{m}^+, 
\boldsymbol{u}_{m}^-,
 \boldsymbol{u}_{m+1}^+, \boldsymbol{u}_{m+1}^-,
 \\
 &\quad ...,\boldsymbol{u}_{d}^+, \boldsymbol{u}_{d}^-
)\ ,
 \end{split}
\]
such that $\boldsymbol{u}_i^{\pm}=\begin{pmatrix}\xi_i\\ \mu_i^{\pm}\xi_i\end{pmatrix}$ for $i=1,2,...,k$ and
$i=k+1,...,l$. Here for $i=1,2,...,k$ we have $\lambda_i<0$, and $\mu_i^{\pm}$ are chosen according to case 1 discussed above;
for $i=k+1,...,l$ we have $0<\lambda_i<\dfrac{\alpha^2}{4}$, and $\mu_i^{\pm}$ are chosen according to case 2 discussed above.
When $i=l+1,...,m$ we have $0<\lambda_i=\dfrac{\alpha^2}{4}$, and in this case $\mu_i^{\pm}=-\dfrac{\alpha}{2}$, so that
$\boldsymbol{u}_i^{+}=\begin{pmatrix}\xi_i\\ -\frac{\alpha}{2}\xi_i\end{pmatrix}$, $\boldsymbol{u}_i^-=\textbf{a}_i$ are chosen according to case 3 discussed above.
When $i=m+1,...,d$ we have $0<\dfrac{\alpha^2}{4}<\lambda_i$ and in this case
$\boldsymbol{u}_i^{\pm}=\begin{pmatrix}\xi_i\\ \mu_i^{\pm}\xi_i\end{pmatrix}$ and $\mu_i^{\pm}$ are chosen according to case 4 discussed above.

Finally by picking the matrix $P$ as above we have
$$P^{-1}AP=\text{diag}(A_1,...,A_d)$$
with each $A_i$ being a $2\times 2$ block matrix, and
$A_i=\begin{pmatrix}\mu_i^+ & 0\\ 0& \mu_i^-\end{pmatrix}$
for $i=1,2,...,k, k+1,...,l, m+1,...,d$;
$A_i=\begin{pmatrix}-\frac{\alpha}{2} & 1\\ 0& -\frac{\alpha}{2}\end{pmatrix}$
for $i=l+1,...,m$.

\end{IEEEproof}

\textbf{Remarks.}

\begin{enumerate}

\item No matter how large $\alpha$ is, the eigenvalues $\mu_1^+, ..., \mu_k^+$ of $A$ are always real and positive.
This means that the number of unstable directions for the Hamiltonian field $\nabla^\perp H$ is always equal to the number of unstable directions
of the Hessian matrix $\nabla^2_X f(X_0)$. In particular, it means that the ``partially damped" friction $b(Y)=-\alpha I^0 Y$
cannot kill all unstable directions of $\nabla^\perp H$, no matter how large $\alpha>0$ is.
\item If $\alpha\rightarrow \infty$, then $\mu_1^+ \rightarrow 0+$.
 It is also very interesting to look at the case of the differential equation for the original Nesterov's method \eqref{Intro:Eq:ODENestrovOriginal}. We see that in this case, the friction coefficient $\alpha$ becomes time--dependent $\alpha(t)= {3}/{t}$. At the beginning, when $t$ is small, we are in the asymptotic regime in which $\alpha>0$ is large, so that the process starts to move according to the unstable direction pointed by the eigenvector that corresponds to $\mu_1^+>0$ in Proposition \ref{prop:LinearAlgebraTechnical}. 
\item In \cite{[ONeill-Wright]}, the authors obtain in their Theorem 4 and Corollary 5 some linear algebra
results that are similar to our Proposition \ref{prop:LinearAlgebraTechnical}.
Yet there are several differences. First, the set--up of \cite{[ONeill-Wright]} does not make use of Hamiltonian structure, and we have revealed
the Hamiltonian structure behind the heavy--ball scheme. Second, we are taking advantage of
the continuous approximation, so that we can make use of very delicate results of Kifer about exit behavior (see \cite{[Kifer1981]} and Theorem \ref{Thm:ExitOneSaddlePoint}).
Third, the analysis of \cite{[ONeill-Wright]} focuses more on showing that the heavy ball method does not converge to
 saddle point, which is originated from \cite{lee2016gradient}. In our case, we follow \cite{hu2019diffusion}, \cite{hu2017fast} and we reduce the problem to the analysis of the dynamics given by stochastic
differential equations with small diffusion term. Our method in a sense is an attempt to investigate the question about exit from saddle points
raised at the end of \cite{[ONeill-Wright]}.

\end{enumerate}

\section{Proof of Theorem \ref{Thm:Kifer1981Strengthened}}
\label{Appendix:ProofOfTheoremKifer1981Strengthened}
%

%

Below we give a proof of Theorem \ref{Thm:Kifer1981Strengthened}.

\begin{IEEEproof}
Let $h$ be the conjugacy mapping in the Linearization Assumption \ref{Assumption:Linearization}.
In this case, by considering $\mathcal{Y}^\varepsilon(t)=h(Y^\varepsilon(t))$ in which $h$ is the linearization homeomorphism that we discussed above, we can transform 
\eqref{Eq:WeakLimitHamiltonianStandardRandomPerturbationEps:Main} into the following equation
\begin{equation}\label{Thm:Kifer1981Strengthened:Eq:WeakLimitHamiltonianStandardRandomPerturbationEpsLinerized}
d\mathcal{Y}^\varepsilon(t)=[(JQ-\alpha I^0) \mathcal{Y}^\varepsilon(t)+ \eps \Psi(\mathcal{Y}^\varepsilon(t))]dt+\sqrt{\eps}\Sigma(\mathcal{Y}^\varepsilon(t))dW_t \ .
\end{equation}
Here the additional term $\Psi$ is smooth in $\mathcal{Y}$. As before, we set $A=JQ-\alpha I^0$ and we have the mild solution of \eqref{Thm:Kifer1981Strengthened:Eq:WeakLimitHamiltonianStandardRandomPerturbationEpsLinerized} written as

\begin{equation*}
\begin{split}
\mathcal{Y}^\varepsilon(t)=e^{At}\mathcal{Y}^\varepsilon(0) + \sqrt{\eps}\int_0^t e^{A(t-s)}\Sigma(\mathcal{Y}^{\varepsilon}(s))dW_s 
\\
+ \eps \int_0^t e^{A(t-s)}\Sigma(\mathcal{Y}^\varepsilon(s))ds \ .
\end{split}
\end{equation*}

Let us pick an orthonormal basis $\zeta_1,...,\zeta_{2d}$ in $\mathbb{R}^{2d}$ such that
$\zeta_i=\dfrac{1}{\sqrt{1+(\mu_i^+)^2}}\begin{pmatrix}\xi_i \\ \mu_i^+ \xi_i \end{pmatrix}$ for $i=1,2,...,k$.
Set the orthogonal matrix
\begin{equation*}
M=\begin{pmatrix} \zeta_1 & \zeta_2 & ... & \zeta_k & \zeta_{k+1} & ... & \zeta_{2d} \end{pmatrix} \ .
\end{equation*}
Then we have
\begin{equation*}
AM=M\begin{pmatrix} \text{diag}(\mu_1^+,...,\mu_k^+) & 0 \\ 0 & \widetilde{A} \end{pmatrix}=M\widehat{A} \ ,
\end{equation*}
in which $\widetilde{A}$ is a matrix of size $(2d-k)\times (2d-k)$, and
$\widehat{A}=\begin{pmatrix} \text{diag}(\mu_1^+,...,\mu_k^+) & 0 \\ 0 & \widetilde{A} \end{pmatrix}$.
In this way we have
\begin{equation*}
e^{At}=M^{T}\begin{pmatrix} \text{diag}(e^{\mu_1^+ t},...,e^{\mu_k^+ t}) & 0 \\ 0 & e^{\widetilde{A}t} \end{pmatrix}M
=M^Te^{\widehat{A} t}M \ ,
\end{equation*}
so that
$$e^{\widehat{A}t}=\begin{pmatrix} \text{diag}(e^{\mu_1^+ t},...,e^{\mu_k^+ t}) & 0 \\ 0 & e^{\widetilde{A}t} \end{pmatrix} \ .$$

Consider $\mathcal{Y}^\varepsilon(t)=\sum\limits_{i=1}^{2d} y_i^\varepsilon(t)\zeta_i$ for some $y^\varepsilon(t)=(y_1^\varepsilon(t),...,y_{2d}^\varepsilon(t))^T$
in $\mathbb{R}^{2d}$ with fixed $\varepsilon>0$ and $t\geq 0$. Then the equation for
$y^\varepsilon(t)$ in terms of mild solution takes the form

\begin{equation*}
\begin{split}
 y^\varepsilon(t)=e^{\widehat{A} t}y^\varepsilon(0)
+ \varepsilon\int_0^t
e^{\widehat{A}(t-s)}\Sigma(y^{\varepsilon}(s))dW_s
\\
+ \varepsilon^2\int_0^t e^{\widehat{A}(t-s)}\Psi(y^\varepsilon(s))ds \ .
\end{split}
\end{equation*}

Write $\Sigma(y)=\begin{pmatrix} \Sigma_1(y) \\ ... \\ \Sigma_k(y) \\ \widetilde{\Sigma}(y)\end{pmatrix}$
in which $\widetilde{\Sigma}(y)$ is a matrix of size $(2d-k)\times 2d$,
$\Psi(y)=\begin{pmatrix}\Psi_1(y) \\ ... \\ \Psi_k(y) \\ \widetilde{\Psi}(y) \end{pmatrix}$
in which $\widetilde{\Psi}(y)$ is a column vector of size $2d-k$, and
$y^\varepsilon(t)=\begin{pmatrix}y^\varepsilon_1(t) \\ ... \\ y^\varepsilon_k(t) \\ \widetilde{y}^\varepsilon(t) \end{pmatrix}$
in which $\widetilde{y}^\varepsilon(t)$ is a column vector of size $2d-k$. Then we have, for $i=1,2,...,k$,

\begin{equation}\label{Thm:Kifer1981Strengthened:Eq:MildSolutionLiearizedEquationCanonicalCoordinate:iFrom1Tok}
\begin{split}
y^\varepsilon_i(t)=e^{\mu_i^+ t}y^\varepsilon_i(0)
+ \sqrt{\eps}\int_0^t
e^{\mu_i^+(t-s)}\Sigma_i(y^{\varepsilon}(s))dW_s
\\
+ \eps \int_0^t e^{\mu_i^+(t-s)}\Psi_i(y^\varepsilon(s))ds \ ,
\end{split}
\end{equation}
and

\begin{equation*} 
\begin{split}
\widetilde{y}^\varepsilon(t)=e^{\widetilde{A} t}\widetilde{y}^\varepsilon(0)
+ \sqrt{\eps} \int_0^t
e^{\widetilde{A}(t-s)}\widetilde{\Sigma}(y^{\varepsilon}(s))dW_s
\\
+ \eps \int_0^t e^{\widetilde{A}(t-s)}\widetilde{\Psi}(y^{\varepsilon}(s))ds \ .
\end{split}
\end{equation*}

By the spectral radius theorem (see \cite{[YosidaFunctionalAnalysis]}) we know that

\begin{equation*}\label{Thm:Kifer1981Strengthened:Eq:SpectralRadiusTheorem}
\lim\limits_{s\rightarrow\infty}\|e^{\widetilde{A}s}\|^{1/s}\leq e^{-\mu^- t}
\end{equation*}
for some $\mu^->0$. This implies that there is some positive constant $C^{(1)}>0$ such that
\begin{equation*}\label{Thm:Kifer1981Strengthened:Eq:Kifer1981:Eq3-9}
\|e^{\widetilde{A}s}\|\leq C^{(1)}e^{-\frac{\mu^-}{2}s}
\end{equation*}
for all $s\geq 0$. From here, we can argue with the same reasoning as in \cite[Section 8]{[BakhtinHeteroclinic]}
and \cite[Appendix A]{hu2017fast}. We see that, to analyze the exit time and exit distribution of the process
$y^\varepsilon(t)$, it suffices to look at the vector consisting of the first $k$--components
$\mathfrak{y}^\varepsilon(t)=(y_1^\varepsilon(t),...,y_k^\varepsilon(t))^T$. The latter process $\mathfrak{y}^\varepsilon(t)$ is governed by the
equations \eqref{Thm:Kifer1981Strengthened:Eq:MildSolutionLiearizedEquationCanonicalCoordinate:iFrom1Tok} where $i=1,2,...,k$.
We can then make full use of the arguments in \cite[Section 8]{[BakhtinHeteroclinic]} and \cite[Appendix A]{hu2017fast}, and we conclude
the statement of this Theorem.
\end{IEEEproof}

\section{Proof of Proposition \ref{prop:LyapunovFunctionFiniteTravelTime}}
\label{Appendix:ProofPropositionLyapunovFunctionFiniteTravelTime}

\begin{IEEEproof}
Let us take the Lyapunov function as the Hamiltonian
$$H(X(t), V(t))=\dfrac{1}{2}(V(t))^2+f(X(t)) \ .$$
Then along the flow of
\eqref{Eq:HeavyBallStrongLimitHamiltonian}
we have

\begin{equation}\label{Eq:Newton-LeibnizFormulaLyapunovFunctionISHamiltonian}
H(X(t), V(t))-H(X(0), V(0))
= -\alpha \int_0^t (V(s))^2 ds \ .
\end{equation}

From the equation \eqref{Eq:Newton-LeibnizFormulaLyapunovFunctionISHamiltonian} we see that when the process $Y(t)$
is a distance away from the $X$--axis, the
Hamiltonian function $H(X(t), V(t))$ is strictly decaying. However, if the process $Y(t)$ crosses the $X$--axis, then $V(s)$ takes
$0$ value along the trajectory and it is not guaranteed that the Hamiltonian $H(X(t), V(t))$ keeps decaying.
It is in this aspect that we see the effect
of interacting the friction with the Hamiltonian flow.

A crossing through $X$--axis happens in the following two cases: (a) It is outside a neighborhood of a critical point $O$
of the Hamiltonian flow $\nabla^\perp H$, that is either a saddle point or a local minimum.
Then the friction vanishes at the crossing point, but the Hamiltonian flow is not zero there, since $|\nabla^\perp H|\geq K>0$,
and thus it is the Hamiltonian flow that brings the process
immediately to a region where $(V(s))^2$ is still strictly positive. Combining this with
\eqref{Eq:Newton-LeibnizFormulaLyapunovFunctionISHamiltonian} we see that even
in this case the Hamiltonian $H(X(t), V(t))$ keeps strictly decaying along the flow of $Y(t)$;
(b) The flow of $Y(t)$ approaches a critical point $O$ of the Hamiltonian flow $\nabla^\perp H$, that is either a saddle point or a local minimum.
In this case, we are entering a neighborhood of the critical point.

Suppose we start the process $Y(t)$ in \eqref{Eq:HeavyBallStrongLimitHamiltonian}
from an initial point $Y_0=(X_0, V_0)$ that stays away from a neighborhood of a saddle point $O$ of the Hamiltonian flow
$\nabla^\perp H$. From the above reasoning we see that, in finite time $T_0>0$ the process $Y(t)$
reaches a neighborhood of another critical point of the Hamiltonian flow $\nabla^\perp H$, that is either a
local minimum point or a saddle. Combining this with the arguments that lead to Lemma 3.3 in \cite{hu2017fast}, we
conclude the validity of this proposition.
\end{IEEEproof}

\section{Proof of Theorem \ref{Thm:MajorResultHamiltonianForm}}
\label{Appendix:ProofTheoremMajorResultHamiltonianForm}

\begin{IEEEproof}
The proof of this Theorem follows the same lines of arguments as those in the proof of Theorem 3.4 in \cite{hu2017fast}.
As $Y_t^\varepsilon$ is a strong Markov process, we see that each of $\sigma_j^\eps -\tau_j^\eps $ in distributed in the same way
as $\tau^\varepsilon_{(x,v)}$ in \eqref{Eq:HittingTimeNeighborhoodSaddle}. However, we note that in this case, the initial condition
$(x,v)$ in $\tau^\varepsilon_{(x,v)}$ will in general be random. In this case, based on
Theorem \ref{Thm:Kifer1981Strengthened}, which is a uniform version of Theorem \ref{Thm:ExitOneSaddlePoint},
one can show that each of $\sigma_j^\eps -\tau_j^\eps$ is distributed in such a way that for any $r>0$
there exist some $\varepsilon_0^{(1)}>0$ such that for any $0<\varepsilon<\varepsilon_0^{(1)}$ we have
\begin{equation}\label{Thm:MajorResultHamiltonianForm:Eq:EstimateTauToSigma}
\dfrac{\mathbf{E}(\sigma_j^{\eps}-\tau_j^{\eps})}{\ln(\varepsilon^{-1})}\leq \dfrac{1}{2\mu_0^{(i)}}+r \ .
\end{equation}
Here $\mu_0^{(i)}=\dfrac{-\alpha+\sqrt{\alpha^2-4\lambda_1^{(i)}}}{2}$.
Notice that by Definition \ref{Def:StrictSaddleProperty}, we have $-\lambda_1^{(i)}\geq \gamma_1$, so that
$\mu_0^{(i)}\geq \dfrac{-\alpha+\sqrt{\alpha^2+4\gamma_1}}{2}$. Thus
\begin{equation}\label{Thm:MajorResultHamiltonianForm:Eq:BoundReciprocalMu0i}
\dfrac{1}{\mu_0^{(i)}}\leq \dfrac{2}{-\alpha+\sqrt{\alpha^2+4\gamma_1}}
=\dfrac{1}{2\gamma_1}\left(\sqrt{\alpha^2+4\gamma_1}+\alpha\right) \ .
\end{equation}

Moreover, it is not difficult to show that there exists some $\varepsilon_0^{(2)}>0$
such that for any $0<\varepsilon<\varepsilon_0^{(2)}$, 
\begin{equation}\label{Thm:MajorResultHamiltonianForm:Eq:EstimateSigmaToTau}
\mathbf{E}(\tau_j^{\eps}-\sigma_{j-1}^{\eps})\leq C
\end{equation}
for some constant $C>0$. Now we decompose
\begin{equation}\label{Thm:MajorResultHamiltonianForm:Eq:ExpansionTotalConvergenceTime:Infinite}
\mathcal{T}_{(X,V)}^{H,\varepsilon}=\sigma_0^{\eps}+(\tau_1^{\eps}-\sigma_0^{\eps})+(\sigma_1^{\eps}-\tau_1^{\eps})+... \ .
\end{equation}
We notice the fact that when the deterministic process \eqref{Eq:HeavyBallStrongLimitHamiltonian} leaves
each of the $G_i$, it never returns to the same $G_i$. Therefore, the expansion
\eqref{Thm:MajorResultHamiltonianForm:Eq:ExpansionTotalConvergenceTime:Infinite}
will terminate after passing through at most $k$--saddle points, i.e,

\begin{equation}\label{Thm:MajorResultHamiltonianForm:Eq:ExpansionTotalConvergenceTime:k}
\mathcal{T}_{(X,V)}^{H,\varepsilon}=\sigma_0^{\eps}+(\tau_1^{\eps}-\sigma_0^{\eps})+(\sigma_1^{\eps}-\tau_1^{\eps})+...+(\tau_k^{\eps}-\sigma_{k-1}^{\eps})+(\sigma_k^{\eps}-\tau_k^{\eps})+(\tau_{k+1}^{\eps}-\sigma_{k}^{\eps}) \ .
\end{equation}
We can then see the validity of this theorem from \eqref{Thm:MajorResultHamiltonianForm:Eq:ExpansionTotalConvergenceTime:k}, \eqref{Thm:MajorResultHamiltonianForm:Eq:EstimateTauToSigma}, \eqref{Thm:MajorResultHamiltonianForm:Eq:EstimateSigmaToTau}
and \eqref{Thm:MajorResultHamiltonianForm:Eq:BoundReciprocalMu0i}.
\end{IEEEproof}

%
%
%

\end{document}